\documentclass[11pt,thmsa]{article}

\usepackage{amssymb,latexsym}
\usepackage{amsmath,amscd}
\usepackage[all]{xy}

 \usepackage[latin1]{inputenc}

 \usepackage{mathrsfs}
 \usepackage{dsfont}
 \usepackage{amsmath}
 \usepackage{amsthm}
 \usepackage{amsfonts}
 \usepackage{amssymb}
 \usepackage[pdftex]{graphicx}
 \usepackage{wrapfig}
 \usepackage{layout}
 \usepackage{verbatim}
 \usepackage{alltt}
 \usepackage{xypic}
 \usepackage{bbm}
 \usepackage{tikz}
 \input xy
 \xyoption{all}

 \DeclareMathAlphabet{\mathpzc}{OT1}{pzc}{m}{it}

 \newtheorem{theorem}{Theorem}[section]
 \newtheorem{lemma}[theorem]{Lemma}
 \newtheorem*{theorem1}{Theorem }
 \newtheorem*{theorem3}{Theorem \ref{theorem:integration_everything}}
 
 \newtheorem{proposition}[theorem]{Proposition}
 
 \newtheorem{corollary}[theorem]{Corollary}
 \newtheorem{definition}[theorem]{Definition}
  \theoremstyle{definition}
 \newtheorem{example}[theorem]{Example}
\newtheorem{examples}[theorem]{Examples}
 \newtheorem{remark}[theorem]{Remark}

\newtheorem*{acknowledgements}{Acknowledgements}
\renewenvironment{proof}{\noindent{\it
Proof.}}{\bgroup\hspace{\stretch{1}}$\square$\egroup\medskip\par}

\newcommand{\Rep}{\textrm{Rep}}

\newcommand{\RRep}{\mathcal{R}\textrm{ep}^{\infty}}
\newcommand{\URRep}{\mathcal{\hat{R}}\textrm{ep}^{\infty}}

\newcommand{\ad}{\mathrm{ad}}

\newcommand{\id}{\mathrm{id}}
\newcommand{\End}{\mathrm{End}}
\newcommand{\Hom}{\mathrm{Hom}}
\newcommand{\RHom}{\underline{\mathrm{Hom}}}

\newcommand{\smooth}{\mathcal{C}^{\infty}}
\newcommand{\Path}{\mathsf{P}}
\newcommand{\ev}{\textrm{ev}}

\newcommand{\Chen}{\mathsf{C}}
\newcommand{\B}{\mathsf{B}}
\newcommand{\s}{\mathsf{s}}
\newcommand{\front}{P}
\newcommand{\back}{Q}
\newcommand{\A}{\mathsf{A}}
\newcommand{\D}{D}
\newcommand{\Hol}{\mathsf{Hol}}
\begin{document}

\def\tetra{\draw(A) -- (B) -- (C) -- (D) -- cycle; \draw[shade](A)--(B)--(D); \draw[shade](B)--(C)--(D);\draw(B)--(D); \draw[dotted](A)--(C);\draw[dotted](A)--(C); \draw(A)--(D);}
\def\tetraA{\draw(A) -- (B) -- (C) -- (D) -- cycle; \draw[shade=gray](B)--(C)--(D);\draw[thin](B)--(D); \draw[dotted](A)--(C);\draw[dotted](A)--(C); \draw[thick](A)--(B);}
\def\tetraB{\draw(A) -- (B) -- (C) -- (D) -- cycle; \shade(A)--(C)--(D);\draw(B)--(D); \draw(A)--(D); \draw[dotted](A)--(C);\draw(C)--(D);}
\def\tetraC{\draw(A) -- (B) -- (C) -- (D) -- cycle; \draw[shade=gray](A)--(B)--(D);\draw(B)--(D); \draw[dotted](A)--(C); \draw[dotted](A)--(C); \draw(A)--(D);\draw[thick](D)--(C);}
\def\tetraD{\draw(A) -- (B) -- (C) -- (D) -- cycle; \draw[shade=gray](A)--(B)--(C);\draw(B)--(D); \draw[dotted](A)--(C);\draw[dotted](A)--(C);\draw(A)--(D);  }

\def\triangle{\draw(A) -- (B) -- (C)--cycle; \draw[shade](A)--(B)--(C);\draw(A)--(C);}
\def\triangleA{\draw(A) -- (B) -- (C)  -- cycle;\draw[ultra thick](A)--(B);\draw[ultra thick](C)--(B);}
 \def\triangleB{\draw(A) -- (B) -- (C)  -- cycle;\draw[ultra thick](A)--(C);}


\def\line{\draw[ultra thick](A) -- (B);\fill[ultra thick](A);}


\def\bleft{\draw(X) -- (Y) -- (Z)--(W);}
\def\bright{\draw(x) -- (y) -- (z)--(w);}

\def\dtetra{\draw(A) -- (B) -- (C) -- (D) -- cycle; \draw(A)--(B)--(D); \draw(B)--(C)--(D);\draw(B)--(D); \draw(A)--(C);\draw[dotted](A)--(C); \draw(A)--(D);}
\def\dtetraA{\draw(A) -- (B) -- (C) -- (D) -- cycle; \draw[shade=gray](B)--(C)--(D);\draw(B)--(D); \draw(A)--(C);\draw(A)--(C); \draw[thick](A)--(B);}
\def\dtetraB{\draw(A) -- (B) -- (C) -- (D) -- cycle; \shade(A)--(B)--(D);\draw[thick](D)--(C); \draw(A)--(D); \draw(A)--(C);\draw(C)--(D);\draw(D)--(B);}
\def\dtetraC{\draw(A) -- (B) -- (C) -- (D) -- cycle; \draw[shade=gray](A)--(C)--(D);\draw(B)--(D); \draw(A)--(C); \draw(A)--(C); \draw(A)--(D);\draw(B)--(C);}
\def\dtetraD{\draw(A) -- (B) -- (C) -- (D) -- cycle; \draw[shade](A)--(B)--(C);\draw(C)--(A); }

\vspace{15cm}
 \title{The $\mathrm{A}_\infty$ de Rham theorem and integration of representations up to homotopy}
\author{Camilo Arias Abad\footnote{Institut f\"ur Mathematik, Universit\"at Z\"urich,
camilo.arias.abad@math.uzh.ch. Partially supported
by SNF Grant 200020-121640/1, the Forschungskredit of the Universit\"at Z\"urich and a Junior Research Fellowship from ESI in Vienna. } \hspace{0cm} and
Florian Sch\"atz\footnote{Center for Mathematical Analysis, Geometry and Dynamical Systems, IST Lisbon, fschaetz@math.ist.utl.pt. 
Partially supported by FCT/POCTI/FEDER through project PTDC/MAT/098936/2008, post-doctoral grant SFRH/BPD/69197/2010 
and a Junior Research Fellowship from ESI in Vienna.}}

 \maketitle
 \begin{abstract} 
We use Chen's iterated integrals to integrate representations up to homotopy.
That is, we  construct an $\mathsf{A}_\infty$ functor $\int : \RRep(A) \rightarrow \URRep(\Pi_{\infty}(A))$ from the representations up to homotopy of a Lie algebroid $A$ to those of
its infinity groupoid. This construction extends the usual integration of representations in Lie theory. We discuss several examples including
Lie algebras and Poisson manifolds.

The construction is based on an $\A_\infty$ version of de Rham's theorem due to Gugenheim \cite{Gugenheim}.
The integration procedure we explain here amounts to extending the construction of parallel transport for superconnections, introduced by Igusa \cite{I} and Block-Smith \cite{BS},  to the case of certain differential graded manifolds.
\end{abstract}

\tableofcontents

\section{Introduction}\label{s:intro}

The purpose of this work is to understand Lie theory for representations up to homotopy. That is, we study the way in which the $\mathrm{L}_\infty$ representations of a Lie algebroid are related to the $\A_\infty$ representation of the corresponding global object. 

Let us first consider the case of ordinary representations. If $G$ is the source simply connected
Lie groupoid integrating the Lie algebroid $A$, then the categories of representations $\Rep(G)$ and $\Rep(A )$ are naturally equivalent. From the point of view of differential graded manifolds, this can be seen as an instance of the fact that flat connections on $\mathcal{M}$ correspond to representations of the fundamental groupoid $\Pi_1(\mathcal{M})$. Indeed, there is a precise sense in which $G$ is the fundamental groupoid of the differential graded manifold $A[1]$ associated to $A$. Hence, ordinary representations only capture information about the 
fundamental groupoid of the associated differential graded manifold. On the other hand, representations up to homotopy capture information about the higher homotopy group(oid)s as well.

Igusa \cite{I} and Block-Smith \cite{BS} showed that flat $\mathbb{Z}$-graded superconnections on a smooth manifold $M$ correspond to representations up to homotopy of its infinity groupoid, namely, the simplicial set of smooth singular chains on $M$.
Building on this, we will show that for any Lie algebroid $A$, there is an $\A_\infty$ functor $\int$ from the category of representations up to homotopy of $A$ to the category of unital representations
up to homotopy of the infinity groupoid of the differential graded manifold $A[1]$ associated to $A$. Moreover, we show that, ultimately, the $\A_\infty$ functor $\int$ relies on an $\A_\infty$ version of
de Rham's theorem that is due to Gugenheim \cite{Gugenheim}, which is constructed using the theory of Chen's iterated integrals.

It has been observed in \cite{AC1,AC2} that some of the natural representations associated to Lie groups and Lie algebras  appear as representations up to homotopy in the case of Lie groupoids and Lie algebroids. These constructions provide many examples to which the $\A_\infty$  functor $\int$ can be applied. We hope to come back to those examples, particularly in the case of Poisson manifolds, in a sequel.

The paper is organized as follows. In the remaining part of the introduction we will  make some comments on simplicial sets, describe in more detail the question that we study and explain in general terms how the $\A_\infty$ functor 
$\int$ is constructed. In Section \S \ref{section:reps_and_simplices} we review Chen's iterated integrals and give a proof of the $\A_\infty$ de Rham 
theorem,  originally due to Gugenheim \cite{Gugenheim}. In Section \S \ref{section:functor} we show that general considerations about $\A_\infty$ algebras
allow one to construct the integration $\A_\infty$ functor $\int$ using the $\A_\infty$ de Rham theorem. In Section \S \ref{section:examples} we discuss some examples to which
the construction applies, including smooth manifolds, Lie algebras and Poisson manifolds. We also show that, when applied to ordinary representations, the $\A_\infty$ functor $\int$ specializes to
the usual integration of representations in Lie theory.
In the Appendix we review the definitions and some general facts concerning representations up to homotopy of Lie algebroids and $\A_\infty$-algebras.

\subsection{General comments on simplicial sets}

It is well known that the theory of simplicial sets allows one to look at some aspects of homotopy theory, higher category theory and Lie theory in a unified way. Let us briefly recall this point of view.
We denote by $\Delta$ the category of nonempty finite ordered sets, and by $[n]$ the ordered set

\[ [n]=\{0\leq 1\leq \dots \leq n\},\] 
seen as a category. We are interested in simplicial sets that appear via the following constructions:

\begin{example}\label{examplenerve}
The nerve functor
\[\mathsf{N}: \mathsf{Cat} \rightarrow \mathsf{Set}_{\Delta}, \]

from the category of small categories to the category of simplicial sets is defined by
 \[\mathsf{N}(\mathcal{C})_n:=\mathrm{ Hom}_{\mathsf{Cat}}([n],\mathcal{C}).\]
In this way, one can view the category $\mathsf{Cat}$ as a full subcategory of the category $\mathsf{Set}_{\Delta}$, and it is easy to characterize precisely the simplicial sets
that are nerves of categories in terms of horn filling conditions.

\end{example}

\begin{example}\label{examplesingular}
We will denote by $\Delta_n$ the geometric $n$-simplex:
\[ \Delta_n:=\{(t_1,...,t_n)\in \mathbb{R}^{n}: 1 \ge t_1 \ge t_2 \ge \cdots \ge t_n \ge 0\} .\]

The sequence of spaces $\Delta_{\bullet}$, together with the usual face and degeneracy maps, forms a cosimplical space.  There is a functor
 $Sing: \mathsf{Top} \rightarrow \mathsf{Set}_{\Delta}$ from topological spaces to simplicial sets given by the formula 
\[Sing(X)_n:=\mathrm{Hom}_{\mathsf{Top}}(\Delta_n, X).\]
The simplicial sets obtained by this procedure are always Kan complexes, also called infinity groupoids (see for example \cite{BV, Jo, Lu1}). This means that they satisfy
some horn filling conditions similar to those that characterize the nerves of categories.

\end{example}

\begin{example}\label{examplesullivan}
Let us also consider the spatial realization functor from rational homotopy theory introduced by Sullivan in \cite{Sul}. This is a functor
\[ \mathsf{S}: \mathsf{DGCA} \rightarrow \mathsf{Set}_{\Delta},\]
from the category of differential graded commutative algebras to simplicial sets.
It is defined by:
\[ \mathsf{S}(\Omega)_n= \mathrm{Hom}_{\mathsf{DGCA}}(\Omega,\Omega(\Delta_n)),\]
where $\Omega(\Delta_n)$ denotes the algebra of differential forms on $\Delta_n$.

\end{example}

The homotopy type of a space $X$ can be reconstructed from the simplicial set $Sing(X)$ and, as we mentioned before,
any category can be reconstructed from its nerve. Thus, one can think of the theory of simplicial sets as a simultaneous generalization
of category theory and homotopy theory. An explanation of the role of simplicial sets in higher category theory can be found in the first chapter of \cite{Lu1}.

Given a space $X$, there is always a map of simplicial sets 
\[\pi: Sing(X) \rightarrow \mathsf{N}(\Pi_1(X)),\]
which sends a simplex to the sequence of homotopy classes of paths associated to its edges. Thus, the simplicial set $Sing(X)$ is a refinement
of the fundamental groupoid, which captures information not only about homotopy classes of paths, but also about the higher homotopy group(oid)s.
Because of this, $Sing(X)$ is often denoted by $\Pi_{\infty}(X)$ and called the infinity 
groupoid of $X$. 

Observe that the relationship between examples \ref{examplesingular} and \ref{examplesullivan} is also very close. Given a smooth manifold $X$, denote by $\mathrm{Sing}(X)$ the simplicial subset of $Sing(X)$, which consists of smooth simplices, and by $\Omega(X)$ the algebra of differential forms on $X$.
Clearly, there is a natural isomorphism of simplicial sets

\[\mathsf{S}(\Omega(X))\cong \mathrm{Sing}(X).\]
 
 Since the inclusion $\mathrm{Sing}(X) \hookrightarrow Sing(X)$ is a homotopy equivalence, one can think of the spatial realization functor $S$ as a generalization
 of the singular chain functor. 

This point of view is quite useful in Lie theory. First, recall that a Lie algebroid structure on a vector bundle $A$ is the same as a differential on $\Omega(A)=\Gamma(\Lambda A^{*})$, which, in turn, is equivalent to equipping the graded manifold $A[1]$ with a cohomological vector field. There are natural identifications
 \[\mathrm{Hom}_{\mathsf{Lie-alg}}(TM, A)\cong \mathrm{Hom}_{\mathsf{dg-man}}(T[1]M, A[1])\cong \mathrm{Hom}_{\mathsf{DGCA}}(\Omega(A), \Omega(M)).\]
 
With these identifications in mind, it becomes natural to think of the simplicial set $ \mathsf{S}(\Omega(A))$ as the space of singular chains in the $dg$-manifold $A[1]$, and to use the notation
\[S(\Omega(A))=\mathrm{Sing}(A[1])=\Pi_{\infty}(A[1]),\]
which, for simplicity, will be denoted by $\Pi_{\infty}(A)$ in the following.

If instead of looking at the infinity groupoid of $A[1]$, one considers the fundamental groupoid $\Pi_1(A[1])$, one obtains a candidate for the groupoid integrating  the Lie algebroid $A$.
Indeed $\Pi_1(A[1])$ is always a topological groupoid
and in case the algebroid is integrable, it is the unique source simply connected Lie groupoid integrating $A$. See \cite{CF,CrF,Severa,Zhu} for more details on this construction.
The integration of $\mathsf{L}_{\infty}$-algebras
can also be interpreted as the computation of higher homotopy groupoids of $dg$-manifolds,
see \cite{G,H,Severa}.

\subsection{The problem}

Our goal is to understand the relation between the global ($\A_\infty$) and the infinitesimal ($\mathsf{L}_\infty$) version of representations up to homotopy.
If $A$ is the Lie algebroid of $G$, one can differentiate any unital representation up to homotopy $E\in \URRep(G)$ to obtain a representation up to homotopy of $A$. This process is explained in
 \cite{AS}, and can be summarized as follows:

\begin{theorem1}
Suppose that $A$ is the Lie algebroid of $G$. Then there is a $dg$-functor 
\[ \Psi:\URRep(G) \rightarrow \RRep(A).\]

\end{theorem1}

The differentiation functor $\Psi$ is constructed by differentiating along the flows of vector fields in all possible directions in order to obtain infinitesimal cochains from global ones.
Appropriately, the inverse process of this iterated differentiation is given by Chen's iterated integrals.

One is led to ask whether it is possible to integrate a representation up to homotopy of $A$ to one of $G$. The answer to this 
question is: no, as explained in Paragraph \ref{subsection:Pi_infty}. Indeed, the global counterpart of a representation up to homotopy of $A$ turns out to be a representation up to homotopy of the infinity groupoid $\Pi_{\infty}(A)$ of $A$, not of its fundamental groupoid $G$.

Observe that a representation up to homotopy of a groupoid $G$ is essentially -- i.e. up to a smoothness condition -- an $\A_\infty$ functor to the $dg$-category of $dg$-vector spaces.
This data can be expressed purely in terms of the simplicial structure of the nerve $\mathsf{N}G$ of $G$, and therefore  it is sensible to speak about representations up to homotopy of 
an arbitrary simplicial set.

Alternatively, one can take the adjoint point of view via Cordier's nerve construction \cite{Co} to define a representation up to homotopy
as a map of simplicial sets, see the Appendix of \cite{BS}. 

Our main result is the following:

\begin{theorem3}
Let $A$ be a Lie algebroid.
The assignments
\begin{align*}
\RRep(A) \ni E  \mapsto \int[E]\in \URRep(\Pi_{\infty}(A))
\end{align*}
and
\begin{eqnarray*}
\int_n: \s \RHom(E_1,E_0) \otimes \cdots \otimes \s \RHom(E_n,E_{n-1}) &\to & \s \RHom(E_n,E_0),\\
\phi_1 \otimes \cdots \otimes \phi_n & \mapsto & \s \left(\Hol(-,\phi_1, \dots,  \phi_n)\right)
\end{eqnarray*}
define an $\mathsf{A}_{\infty}$ functor 
\begin{align*}
\int: \RRep(A) \to \URRep(\Pi_{\infty}(A))
\end{align*}
between the $dg$-category of representations
up to homotopy of $A$ and the $dg$-category of unital representations up to homotopy of $\Pi_{\infty}(A)$.
\end{theorem3}

The construction of the $\A_\infty$ functor $\int$ is based on and inspired by the parallel transport of superconnections studied
by Igusa \cite{I} and Block-Smith \cite{BS}. It seems to us that -- beyond the extension to arbitrary Lie algebroids -- our contributions are as follows:

We show that the parallel transport from \cite{I,BS} can be derived from the $\A_\infty$ version of de Rham's theorem due to Gugenheim \cite{Gugenheim}; this allows us to extend parallel transport to 
an $\A_\infty$ functor.

Moreover, we establish that  the integration process lands in the category of {\em unital} representations up to homotopy. This is significant because, as explained in Paragraph \ref{Some examples}, unitality implies higher dimensional versions of the glueing rules that express how the holonomies behave under a refinement of the triangulation. Observe that the unitality condition also comes up naturally in the construction of the  differentiation functor $\Psi:\URRep(G) \rightarrow \RRep(A)$ from \cite{AS}.

\begin{acknowledgements} Conversations with several people have been of great help during the time we worked on this project. 
We would like to mention Alberto Cattaneo,  Marius Crainic,  Ezra Getzler, Sergei Merkulov, Pavel Mn\"ev, Dmitry Roytenberg, Urs Schreiber, Domingo Toledo  and Bernardo Uribe. We thank Calder Daenzer for mentioning the paper \cite{BS}, and James Stasheff for his valuable comments after reading a preliminary version.
We also thank the Instituto
Superior Tecnico in Lisbon, IMPA in Rio de Janeiro, ESI in Vienna and the Universit\"at Z\"urich for their hospitality during the process of writing this
paper.
\end{acknowledgements}

\section{Simplicial sets and representations up to homotopy}\label{section:reps_and_simplices}

In this section we give the definition of a representation up to homotopy of a simplicial set and discuss some
examples. Later on we will be interested in the case where the simplicial set is the  $\infty$-groupoid $\Pi_{\infty}(A)$ of a Lie algebroid $A$.

\subsection{Definitions}\label{subsection:definitions}

\begin{remark}
Let $X_{\bullet}$ be a simplicial set with face and degeneracy maps denoted by
\[d_i: X_k\rightarrow X_{k-1} \quad \text{and} \quad s_i: X_k \rightarrow X_{k+1},\] 
respectively.  We will use the notation
\begin{eqnarray*}
\front_i&:=&(d_0)^{k-i}:X_k \rightarrow X_i,\\
\back_i&:=&d_{i+1} \circ \dots \circ d_{k} :X_k \rightarrow X_i,
\end{eqnarray*}
for the maps that send a simplex to its $i$-th back and front face.
The $i$-th vertex of a simplex $\sigma \in X_k$ 
will be denoted $v_i(\sigma)$, or simply $v_i$, when no confusion can arise.
In terms of the above operations, one can write
\[v_i=(\front_0\circ \back_{i})(\sigma). \]

A cochain $F$ of degree $k$ on $X_{\bullet}$ with values in an algebra $\A$ is a map
\begin{align*}
F: X_k \to \A.
\end{align*}
As usual, the cup product of two cochains  $F, F'$ of degree $i$ and $j$, respectively, is the cochain of degree $i+j$ defined by the formula:
\[(F \cup F') (\sigma):= F(\back_i(\sigma)) F'(\front_j(\sigma)).\]
\end{remark}

\begin{definition}
A representation up to homotopy of $X_{\bullet}$ consists of the following data:
\begin{enumerate}
\item A graded vector space $E_x=\bigoplus_{k\in \mathbb{Z}}E_x^k$, 
for each zero simplex $x \in X_0$.
\item A sequence of operators $\{F_k\}_{k\geq 0}$, where $F_k$ is a $k$-cochain that assigns to $\sigma \in X_k$ a linear map
\[F_k(\sigma)\in \mathrm{Hom}^{1-k}(E_{v_k(\sigma)},E_{v_0(\sigma)}).\]
\end{enumerate}
These operators are required to satisfy the following  equations:
\begin{equation}\label{structure equations}
\sum_{j=1}^{k-1} (-1)^{j} F_{k-1}(d_j (\sigma))  +
\sum_{j=0}^{k} (-1)^{j+1} (F_{j}\cup F_{k-j})(\sigma)=0. 
\end{equation}
A representation up to homotopy is called unital if
\[ F_1(s_0(\sigma))=\id  \quad \text{and} \quad  F_k(s_i(\sigma))=0 \text{ for } k>1. \]
\end{definition}

\begin{remark}
The definition of representation up to homotopy of a simplicial set is such that when applied to the nerve of a category, it is the same as an $\mathsf{A}_\infty$-functor to the $dg$-category of differential graded vector spaces. These equations were first considered by Sugawara \cite{Sug}. The structure of a representation up to homotopy can also be described from the adjoint point of view via Cordier's nerve construction \cite{Co,Lu1}, as explained in the appendix of \cite{BS}. The case where the simplicial set is the nerve of a Lie groupoid was considered in \cite{AC2}.

Observe that the structure equations of a representation up to homotopy imply that $F_0$ gives each of the vector spaces $E_x$ the structure of a cochain complex.
We will often write $\partial$ instead of $F_0$ for this coboundary operator.
\end{remark}

\begin{remark}
The representations up to homotopy of $X_{\bullet}$ form a $dg$-category: Let $E, E'$ be
two representations up to homotopy of $X_{\bullet}$.
A degree $n$ morphism $\phi \in \RHom^n(E,E')$ is a formal sum
\begin{equation*}
\phi =\phi_0 +\phi_1+ \phi_2 + \dots,
\end{equation*}
where $\phi_k$ is a $k$-cochain which assigns to a simplex $\sigma \in X_k$ a linear map
\[\phi_k(\sigma)\in \mathrm{Hom}^{n-k}(E_{v_k(\sigma)},E'_{v_0(\sigma)}).\]
The differential
\[\D:\RHom^n(E,E')\rightarrow \RHom^{n+1}(E,E')\]
is defined by
\[\D(\phi)=\D(\phi)_0+\D(\phi)_1+ \D(\phi)_2+\dots,\]
where
\begin{equation}\label{differentialonmorphisms}
\D(\phi)_k(\sigma):=\sum_{i+j=k}(-1)^{jn}F_{j}^{'}\cup
 \phi_{i}(\sigma)+ \sum_{i+j=k}(-1)^{n+j+1}\phi_j \cup
F_{i}(\sigma)+\sum_{j=1}^{k-1}(-1)^{j+n}\phi_{k-1} (d_j(\sigma)).
\end{equation}

If $\phi': E' \rightarrow E''$ is a morphism of degree $m$, the composition $ \phi' \circ \phi \in \RHom^{m+n}(E,E'')$ is given by
\[\phi'\circ \phi=(\phi'\circ \phi)_0+ (\phi'\circ \phi)_1+( \phi'\circ \phi)_2 + \dots,\]
where
 \[(\phi'\circ \phi)_k:=\sum_{i+j=k} (-1)^{jn}(\phi_j' \cup \phi_i). \]
 
One easily checks that these operators define a $dg$-category -- the category of representations
 up to homotopy of $X_{\bullet}$ -- which we denote by $\RRep(X_\bullet )$. 

We will be particularly interested in the category of unital representations up to homotopy, 
which we denote by $\URRep(X_\bullet)$. This is the sub-$dg$-category of $\RRep(X_{\bullet})$ whose objects
are unital representations up to homotopy and whose morphisms are unital morphisms, that is, morphisms $\phi$ such that $\phi_k(s_i(\sigma))=0$. 
\end{remark}

\subsection{Some examples}\label{Some examples}

\begin{example}
Let $G$ be a group and $\mathsf{N}G$ the nerve of $G$. A unital representations up to homotopy of $\mathsf{N}G$ on a graded vector space $V$, which is concentrated in degree zero, is the same
as a representation of $G$.
\end{example}
\begin{example}
Let $C$ be a category and $\mathsf{N}C$ the nerve of $C$. A representations up to homotopy of $\mathsf{N}C$ is the same as an $\mathsf{A}_\infty$-functor from the category $\mathbb{R}C$ to the $dg$-category of
complexes of vector spaces. Here $\mathbb{R}C$ denotes the $dg$-category whose objects are those of $C$ and whose morphisms are the linear spans of morphisms of $C$.
\end{example}

\begin{example}
The most important simplicial set for our purposes is the $\infty$-groupoid $\Pi_{\infty}(A)$ of a Lie algebroid $A$, given by
\begin{align*}
(\Pi_{\infty}(A))_k := \textrm{Hom}_{\mathsf{DGCA}}(\Omega(A),\Omega(T\Delta_k)),
\end{align*}
i.e. the morphisms of differential graded commutative algebras between $\Omega(A)$
and the de Rham algebra of the $k$-simplex $\Delta_k$. Observe that this is the same as the set of  Lie algebroid morphisms from $T\Delta_k$ to $A$.
The structure operations for a representation up to homotopy of $\Pi_{\infty}(A)$ are rules that assign
holonomies to simplices, and the structure equations are compatibility conditions between the holonomies. In the pictures below, a shaded face of a simplex denotes the holonomy assigned to it.
For the $1$-simplex we obtain
 
\begin{center}
\begin{tikzpicture}[scale=1]
    \coordinate (A) at (-0.5,0);
    \coordinate (B) at (0.5,0);
    \coordinate (X) at (0,0.5);
    \coordinate (Y) at (-0.1,0.5);
    \coordinate (Z) at (-0.1,-0.5);
    \coordinate (W) at (0,-0.5);
     \coordinate (x) at (-0.5,0.5);
    \coordinate (y) at (-0.4,0.5);
    \coordinate (z) at (-0.4,-0.5);
     \coordinate (w) at (-0.5,-0.5);

\matrix[column sep=0.8cm,row sep=0.5cm]
{
\bleft &node{$\partial ,$} & \line &\bright & node{$=$}&&node{$0.$}\\
};
\end{tikzpicture}
\end{center}

This says that the holonomy assigned to a path should be a chain map between the chain complexes
associated to the endpoints of the paths.

The compatibility for condition for the triangle is\\

\begin{center}
\begin{tikzpicture}[scale=1]
    \coordinate (A) at (-1,-0.5);
    \coordinate (B) at (0,-0.5);
    \coordinate (C) at (0,0.5);
    \coordinate (X) at (0,0.5);
    \coordinate (Y) at (-0.1,0.5);
    \coordinate (Z) at (-0.1,-0.5);
    \coordinate (W) at (0,-0.5);
      \coordinate (x) at (-0.5,0.5);
    \coordinate (y) at (-0.4,0.5);
    \coordinate (z) at (-0.4,-0.5);
     \coordinate (w) at (-0.5,-0.5);

     \matrix[column sep=0.8cm,row sep=0.5cm]
{
\bleft & node{$\partial ,$} &\triangle&\bright & node{$=$}&\triangleA &node[$-$]& \triangleB&  node[$,$]\\
};
\end{tikzpicture}
\end{center}

requiring that -- even though two homotopic paths are assigned different holonomies in general -- any specific homotopy between the two paths induces a homotopy between the corresponding chain maps.

For the tetrahedron one obtains\\

\begin{center}
\begin{tikzpicture}[scale=1]
    \coordinate (A) at (-0.8,-0.3);
    \coordinate (B) at (0,-0.5);
    \coordinate (D) at (-0.1,0.5);
    \coordinate (C) at (0.5,-0.1);
     \coordinate (X) at (0,0.5);
    \coordinate (Y) at (-0.1,0.5);
    \coordinate (Z) at (-0.1,-0.5);
    \coordinate (W) at (0,-0.5);
     \coordinate (x) at (-0.5,0.5);
    \coordinate (y) at (-0.4,0.5);
    \coordinate (z) at (-0.4,-0.5);
     \coordinate (w) at (-0.5,-0.5);

\matrix[column sep=0.8cm,row sep=0.5cm]
{
\bleft & node{$\partial ,$} &\tetra &\bright & node{$=$}&\tetraA &node[$-$]& \tetraC & node{$+$}& \tetraD& node{$-$}& \tetraB &  node[$.$]\\
};
\end{tikzpicture}
\end{center}

\end{example}

\begin{example}
The unitality condition on a representation up to homotopy describes the way in which the holonomies change under a refinement of the triangulation.
For a degenerate one simplex, the unitality implies that the constant path gets assigned the identity. For a degenerate triangle, one obtains
that dividing an interval in two pieces and composing their holonomies gives the holonomy associated to the original interval. 
For a tetrahedron, unitality yields the higher dimensional analogue of this glueing rule.
We illustrate the degenerate tetrahedron as a triangle with a vertex in the middle, which comes from the degeneracy map. Then, the structure equation gives\\

\begin{center}
\begin{tikzpicture}[scale=1]
    \coordinate (A) at (-0.5,-0.5);
    \coordinate (B) at (0.5,-0.5);
    \coordinate (D) at (0,0);
    \coordinate (C) at (0,0.7);
     \coordinate (X) at (0,0.5);
    \coordinate (Y) at (-0.1,0.5);
    \coordinate (Z) at (-0.1,-0.5);
    \coordinate (W) at (0,-0.5);
     \coordinate (x) at (-0.5,0.5);
    \coordinate (y) at (-0.4,0.5);
    \coordinate (z) at (-0.4,-0.5);
     \coordinate (w) at (-0.5,-0.5);

\matrix[column sep=0.8cm,row sep=0.5cm]
{
\dtetraA &node[$-$]& \dtetraB & node{$+$}& \dtetraD& node{$-$}& \dtetraC & node{$=$}&   \bleft & node{$\partial ,$} &\dtetra &\bright node{$=0,$}\\
};
\end{tikzpicture}
\end{center}

where the right hand side vanishes because the tetrahedron under consideration is degenerate. Thus, we obtain\\

\begin{center}
\begin{tikzpicture}[scale=1]
    \coordinate (A) at (-0.5,-0.5);
    \coordinate (B) at (0.5,-0.5);
    \coordinate (D) at (0,0);
    \coordinate (C) at (0,0.7);
     \coordinate (X) at (0,0.5);
    \coordinate (Y) at (-0.1,0.5);
    \coordinate (Z) at (-0.1,-0.5);
    \coordinate (W) at (0,-0.5);
     \coordinate (x) at (-0.5,0.5);
    \coordinate (y) at (-0.4,0.5);
    \coordinate (z) at (-0.4,-0.5);
     \coordinate (w) at (-0.5,-0.5);

\matrix[column sep=0.8cm,row sep=0.5cm]
{
\dtetraD &node[$=$]& \dtetraC & node{$-$}& \dtetraA& node{$+$}& \dtetraB & node[$.$]    \\
};
\end{tikzpicture}
\end{center}

This relation expresses the holonomy of the big triangle in terms of the holonomies associated to the small triangles and
the edges.

\end{example}

\section{Differential forms and singular cochains}\label{section:diff_forms_and_cochains}

In this section, we explain how to use Chen's iterated integrals to construct a natural $\mathsf{A}_\infty$ quasi-isomorphism
\begin{align*}
\psi: (\Omega(M),-d,\wedge) \to (C(M),\delta,\cup),
\end{align*}

between the algebra of differential forms and the algebra of smooth singular cochains of a manifold $M$. 
This result is originally due to Gugenheim  (\cite{Gugenheim}).
Implicitly, this quasi-isomorphism is also present in \cite{BS} and \cite{I}.

\begin{remark}
Given a differential graded algebra $(A,d,\wedge)$, the bar complex
\begin{align*}
\B(\s A) := \bigoplus_{k\ge 1} \left(\s A \right)^{\otimes k}
\end{align*}
carries a coboundary operator $D$ given by
\begin{eqnarray*}
D(\s a_1 \otimes \cdots \s a_n) &:=& \sum_{i=1}^{n} (-1)^{[a_1]+\cdots + [a_{i-1}]} \s a_1 \otimes \cdots \otimes \s a_{i-1} \otimes \s (da_i) \otimes \s a_{i+1} \otimes \cdots \otimes \s a_n\\
&& + \sum_{i=1}^{n-1}(-1)^{[a_1]+ \cdots + [a_i]} \s a_1 \otimes \cdots \otimes \s a_{i-1} \otimes \s (a_i \wedge a_{i+1}) \otimes \s a_{i+1} \otimes \cdots \otimes \s a_n.
\end{eqnarray*}
\end{remark}

The morphism $\psi$ is constructed out of two maps
\begin{align*}
\xymatrix{
\B(\s \Omega(M)) \ar[rr]^{\Chen} && \Omega(\Path M) \ar[rr]^{\mathsf{S}}&& \s C(M).
}
\end{align*}
Here,
\begin{align*}
\xymatrix{
\B(\s \Omega(M)) \ar[rr]^{\Chen} && \Omega(\Path M)
}
\end{align*}
is given by Chen's iterated integrals, while
\begin{align*}
\xymatrix{
\Omega(\Path M) \ar[rr]^{\mathsf{S}}&& \s C(M)
}
\end{align*}
is constructed with the help of a certain family of maps
\begin{align*}
(\Theta_{(k)}: I^{k-1} \to \Path \Delta_k)_{k\ge 1}
\end{align*}
from the cubes to the path spaces of the simplices. The crucial property of this family
is that it relates the cellular structures of the cubes to those of the simplices.
Such a family was already considered by Chen \cite{C}. We will
make use of the explicit family constructed by Igusa in \cite{I}.

\subsection{Chen's iterated integrals}\label{subsection:Chen}

Here, we explain how Chen's iterated integrals provide a map
$$ \Chen: \B(\s \Omega(M)) \rightarrow \Omega(\Path M),$$
from the bar complex of the suspension of $\Omega(M)$ to the differential forms on the path
space of $M$. The most important property of $\Chen$ is that it is almost a chain map
between $\B(\s \Omega(M))$ -- equipped with the bar differential $\overline{D}$ corresponding to the $dg$-algebra $(\Omega(M),-d,\wedge)$ -- and
$\Omega(\Path M)$ -- equipped with the de Rham differential. In fact, $\Chen$ fails to be
a chain map only because of two boundary terms, see Theorem \ref{theorem:Chen}.

\begin{remark}
We use the conventions from \cite{Getzleretal}, except for the definition of the simplices $\Delta_k$, where we follow Igusa (\cite{I}).
\end{remark}

Our first task is to introduce the path space $\Path M$  of $M$, as well as differential forms on it.

\begin{remark}
The path space $\Path M$ of a smooth manifold $M$ is the topological space $\smooth(I,M)$, equipped with the $\mathcal{C}^1$-topology.
By definition, a map $f: X \to \Path M$ is smooth if 
\begin{align*}
f^{\ev}: I \times X \to M, \quad f^{\ev}(t,x):=(f(x))(t) 
\end{align*}
is smooth.

Differential forms on $\Path M$ are defined as follows:
\begin{itemize}
\item[i)] Denote by $\smooth(-,\Path M)$ the category whose objects are pairs $(X,f)$, with $X$ a smooth manifold and $f$ a smooth map
from $X$ to $\Path M$; the morphisms from $(X,f)$ to $(Y,g)$ are all smooth maps $h: X \to Y$ such that $f=g \circ h$.
\item[ii)] Let $\underline{\mathbb{R}}(-)$ be the functor from $\smooth(-,\Path M)$ to the category of real vector spaces
$\mathsf{Vect}$, which maps any object in $\smooth(-,\Path M)$ to $\mathbb{R}$ and every morphism to the identity.
\item[iii)] The functor $\Omega(-)$ from $\smooth(-,\Path M)$ to $\mathsf{Vect}$ is defined via $(X,f) \mapsto \Omega(X)$ and $g \mapsto g^*$.
\item[iv)] A differential form on $\Path M$ is a natural transformation from $\underline{\mathbb{R}}(-)$ to $\Omega(-)$. 
\end{itemize}

Put another way: a differential form $\alpha$ on $\Path M$ is a natural association of a differential form on
$X$, which we will denote by $f^*\alpha$, to a pair $(X,f: X \to \Path M)$.
The idea  is simply that it suffices to know
all the pull backs of a differential form to finite dimensional manifolds in order to know the differential form itself.
The wedge product and the de Rham differential are defined by requiring the pullback operation to preserve them.


We call a continuous map $h: \Path M \to \Path N$ smooth if the composition of any smooth map $f:X\ \to \Path M$ 
with $h$ is smooth. Smooth maps induce pull back maps $\Omega(\Path N) \to \Omega(\Path M)$.
Moreover, observe that every smooth map $h: M \to N$ induces a smooth map
$\Path h: \Path M \to \Path N$ on the path spaces via composition.
\end{remark}

\begin{remark}
Let $X$ be a manifold and $\Delta_k$ the $k$-simplex, i.e.
\begin{align*}
\Delta_k := \{ (t_1,\dots, t_k) \in \mathbb{R}^{k}: 1\ge t_1 \ge t_2 \ge \cdots \ge t_k \ge 0\} \subset \mathbb{R}^{k}.
\end{align*}
We define the push forward along the projection
\begin{align*}
\pi: \Delta_k \times X \to X
\end{align*}
to be the linear map
\begin{align*}
\pi_*: \Omega(\Delta_k \times X) \to \Omega(X)
\end{align*}
of degree $-k$ determined by setting
\begin{align*}
\pi_{*}\left(f(t,x)dt^{i_1}\cdots dt^{i_r}dx^{j_1}\cdots dx^{j_s}\right):=
\left(\int_{\Delta_k} f(t,x) dt^{i_1}\cdots dt^{i_r}\right)dx^{j_1}\cdots dx^{j_s}.
\end{align*}

Observe that for $X$ compact and oriented,
\begin{align*}
\int_{X}\pi_{*}(\alpha) = \int_{\Delta_k \times X}\alpha
\end{align*}
holds for all $\alpha \in \Omega(\Delta_k \times X)$.
\end{remark}

\begin{lemma}\label{lemma:push_forward}
The push forward $\pi_*: \Omega(\Delta_k\times X) \to \Omega(X)$ is a morphism of left $\Omega(X)$-modules of degree $-k$,
i.e. for every $\alpha \in \Omega(X)$ and every $\beta \in \Omega(\Delta_k\times X)$, we have
\begin{align*}
\pi_*(\pi^{*}\alpha \wedge \beta) = (-1)^{|\alpha|k}\alpha \wedge \pi_*(\beta).
\end{align*}

Furthermore, let $\partial \pi$ be the composition
\begin{align*}
\xymatrix{
 \partial \Delta_k \times X \ar[rr]^(0.5){\iota \times \id} && \Delta_k \times X \ar[rr]^(0.5){\pi}&&  X,
}
\end{align*}
where $\partial \Delta_k$ denotes the disjoint union of the codimension one strata
of $\Delta_k$.

Then the following formula holds:
\begin{align*}
\pi_{*} \circ d - (-1)^{k}d \circ \pi_{*} = (\partial \pi)_{*}\circ (\iota \times \id)^{*}.
\end{align*}
Here, the push forward along $\partial \pi$ is understood as the sum over the push forwards
of the individual faces which appear as the connected components of $\partial \Delta_k$.
\end{lemma}

\begin{remark}
Next we define Chen's iterated integrals for differential forms on a smooth manifold $M$.
This is a degree zero linear map:
\begin{align*}
\Chen: \B(\s \Omega(M)):=\bigoplus_{k\ge 1} \left(\s \Omega(M) \right)^{\otimes k} \to \Omega(\Path M).
\end{align*}

Let $\s a_1\otimes \cdots \otimes \s a_n$ be an element of the bar complex $\B(\s \Omega(M))$.
Given any smooth map $f:  I \times X \to M$, we define a differential form on $X$ by the following procedure:
\begin{itemize}
\item[i)] Thinking of $\Delta_n$ as the compactified configuration space of (unordered) points on the interval $I$,
one obtains an extension of $f:  I \times X \to M$ to 
\begin{align*}
f_{(n)}: \Delta_n \times X \to M\times \cdots \times M, \quad f_{(n)}(t_1,\dots, t_n,x):=(f(t_1,x),\dots,f(t_n,x)).
\end{align*}
We will sometimes omit the subscript $(n)$ if the dimension of the simplex is clear from the context.
\item[ii)] Pull back the differential form $a_i$ to $M \times \cdots \times M$ via the $i$'th projection map
\begin{align*}
p_i: M^{\times n} \to M
\end{align*}
and pull back the
wedge product $p_1^{*}a_1\wedge \cdots \wedge p_n^{*}a_n$ to $\Delta_n \times X$ via $f_{(n)}$.
\item[iii)] To obtain a differential form on $X$, push the differential form
$f_{(n)}^{*}(p_1^{*}a_1 \wedge \cdots \wedge p_n^{*}a_n)$ forward to $X$
along
$\pi: \Delta_n \times X \to X$.
\item[iv)] Finally, multiply this differential form with the sign 
\begin{align*}
(-1)^{\sum_{i=1}^n[a_i](n-i)},
\end{align*}
where $[a]$ denotes the degree of an element $\s a \in \s \Omega(M)$.

Note that this is exactly the sign by which $(\s^{-1} \otimes \cdots \otimes s^{-1})(\s a_1 \otimes \cdots \otimes \s a_n)$ differs from $(a_1 \otimes \cdots \otimes a_n)$.
\end{itemize}
We denote the resulting differential form on $X$ by $f^{*}(\Chen(\s a_1\otimes \cdots \otimes \s a_k))$. It is straight forward to check that
this construction is natural and therefore defines
a differential form on the path space $\Path M$.
Let us write down the whole construction more compactly:
\end{remark}

\begin{definition}
Chen's map
\begin{align*}
\Chen: \B(\s \Omega(M)) \to \Omega(\Path M)
\end{align*}
from the bar complex of the suspension of differential forms to the de Rham complex 
of the path space
$\Path M$ is defined by setting
\begin{align*}
f^{*}(\Chen(\s a_1\otimes \cdots \s a_n)):= (-1)^{\sum_{i=1}^n[a_i](n-i)}\pi_{*}((f^{\ev}_{(n)})^{*}(p_1^*a_1\wedge \cdots \wedge p_n^*a_n)) 
\end{align*}
for any object $(X,f)$ in the category $\smooth(-,\Path M)$.
\end{definition}

\begin{remark}

An important property of Chen's map is that the image of an element
$\s a_1\otimes \cdots \otimes \s a_n$ of $\B(\s \Omega(M))$,
for which one of the factors $a_i$ is a function, is zero.
This follows from the observation that in this case the pull back $f^{*}\Chen(\s a_1\otimes \cdots \otimes \s a_n)$ by any smooth map $f: X \to \Path M$ vanishes. Indeed, the differential form
\begin{align*}
(f ^{\ev}_{(n)})^{*}(p_1^{*}a_1\wedge \cdots \wedge p_n^{*}a_n) \in \Omega(\Delta_n \times X)
\end{align*}
is annihilated by the vector field $\frac{\partial}{\partial t^{i}}$, where $t^{i}$ denotes
the $i$-th coordinate on the simplex $\Delta_n$. Hence the push forward of the differential form along $\pi: \Delta_n \times X \to X$ vanishes.
\end{remark}

\begin{lemma}
Chen's map $\Chen$ is natural, i.e. for any smooth map $h: M \to N$ the diagram
\begin{align*}
\xymatrix{
\B(\s \Omega(M)) \ar[r]^{\Chen} & \Omega(\Path M) \\
\B(\s \Omega(N)) \ar[r]^{\Chen} \ar[u]^{\B h} & \Omega(\Path N) \ar[u]_{(\Path h)^*}
} 
\end{align*}
commutes. Here, $\B h$ denotes the linear map which extends
\begin{align*}
\s a_1\otimes \cdots \otimes \s a_n \mapsto \s h^*(a_1)\otimes \cdots \otimes \s h^*(a_n). 
\end{align*}
\end{lemma}

\begin{proof}
This is essentially an exercise in unraveling the definitions. The key-observation is that for any smooth map $f: X\to \Path M$,
\begin{align*}
(\Path h \circ f)^{\ev}_{(n)} = (h\times \cdots \times h) \circ f^{\ev}_{(n)}
\end{align*}
is satisfied.
\end{proof}

\begin{theorem}[Chen]\label{theorem:Chen}
For any element $\s a_1 \otimes \cdots \otimes \s a_n$ of the bar complex $\B(\s \Omega(M))$, the following
identity holds:
\begin{eqnarray*}
d (\Chen(\s a_1\otimes \cdots \otimes  \s a_n)) &=& \Chen(\overline{D}(\s a_1\otimes \cdots \otimes \s a_n)) + \\
&& \hspace{-1cm} ev_1^{*}(a_1) \wedge \Chen(\s a_2 \otimes \cdots \otimes \s a_n) - (-1)^{[a_1]+\cdots + [a_{n-1}]} \Chen(\s a_1\otimes \cdots \otimes \s a_{n-1}) \wedge \ev_0^{*}(a_n).
\end{eqnarray*}
Here:
\begin{enumerate}
\item $\overline{D}$ is the differential corresponding to the differential graded algebra
$(\Omega(M),-d,\wedge)$.
\item  The maps $\ev_0$ and $ev_1:\Path M\rightarrow M$ are given by $\gamma \mapsto \gamma(i)$ for $i=0$ or $1$, respectively.
\end{enumerate}
\end{theorem}

\begin{proof}
It suffices to establish the formula for the pull back of $\Chen(\s a_1\otimes \cdots \otimes \s a_{n})$
by any smooth map $f: X \to \Path M$. 
Using Lemma \ref{lemma:push_forward}, one sees that
$f^{*}d \Chen(\s a_1 \otimes \cdots \otimes \s a_n)$ is equal to
\begin{eqnarray*}
&&(-1)^{n}(\pi_*d\left(f_{(n)}^{*}(p_1^{*}a_1\wedge \cdots \wedge p_n^{*}a_n)\right)) +\\
&&(-1)^{n+1}((\partial \pi)_*(\iota \times \id)^{*}f_{(n)}^{*}(p_1^{*}a_1\wedge \cdots
\wedge p_n^{*}a_n))
\end{eqnarray*}
times the sign $(-1)^{\sum_{i=1}^{n}[a_i](n-i)}$.
The first term gives
\begin{align*}
\sum_{i=1}^{n}(-1)^{[a_1]+\cdots+[a_{i-1}]} f^{*}\Chen(\s a_1 \otimes \cdots \otimes \s(-d a_i) \otimes \cdots \otimes  \s a_n),
\end{align*}
while the second one yields
\begin{eqnarray*}
&& \hspace{-0.5cm} \sum_{i=1}^{n-1}(-1)^{[a_1]+\cdots +[a_i]}f^{*}\Chen(\s a_1\otimes \cdots \otimes \s(a_i\wedge a_{i+1}) \otimes \cdots \otimes \s a_n)\\
&& + (\ev_1\circ f)^{*}a_1 \wedge f^{*}\Chen(\s a_2\otimes \cdots \otimes \s a_n)\\
&& - (-1)^{[a_1]+\cdots + [a_{n-1}]} f^{*}\Chen(\s a_1\otimes \cdots \otimes \s a_{n-1}) \wedge (\ev_0 \circ f)^{*}a_n.
\end{eqnarray*}
\end{proof}

\begin{definition}
The set
\begin{align*}
\smooth_+(I,\partial I) := \{\phi: I \to I \textrm{ smooth and monotone}: \phi(0)=0, \phi(1)=1\} 
\end{align*}
is a monoid under composition and acts on $\smooth(X,\Path M)$ via reparametrizations, i.e. via
\begin{eqnarray*}
\smooth_+(I,\partial I) \times \smooth(X,\Path M) &\to& \smooth(X,\Path M),\\
(\phi,f) &\mapsto& (\phi\bullet f)^{\ev}(t,x):=f^{\ev}(\phi(t),x).
\end{eqnarray*}

A differential form $\alpha \in \Omega(\Path M)$ is reparametrization invariant if for any smooth map
$f: X \to \Path M$ and any reparametrization $\phi \in \smooth_+(I,\partial I)$, the equation
\begin{align*}
f^*\alpha = (\phi\bullet f)^* \alpha 
\end{align*}
holds.
We denote the subcomplex of invariant differential forms by $\Omega_{\textrm{inv}}(\Path M)$.
\end{definition}

\begin{lemma}
The image of Chen's map
\begin{align*}
\Chen: \B(\s \Omega(M)) \to \Omega(\Path M)
\end{align*}
is contained in the subcomplex $\Omega_{\textrm{inv}}(\Path M)$ of reparametrization invariant differential forms on $\Path M$.
\end{lemma}

\begin{proof}
Pick $f: X\to \Path M$ and $\phi \in \smooth_+(I,\partial I)$ arbitrary. We have to check that
\begin{align*}
f^*\Chen(\s a_1\otimes \cdots \otimes \s a_n) = (\phi\bullet f)^*\Chen(\s a_1\otimes \cdots \otimes \s a_n)
\end{align*}
holds for any $a_1,\dots,a_n \in \Omega(M)$.

Using
\begin{align*}
(\phi\bullet f)^{\ev}_{(n)} = f^{\ev}_{(n)} \circ ( \phi \times \cdots \times \phi \times \id): \Delta_n \times X \to M^{\times n},
\end{align*}
one sees that it suffices to prove the equation
\begin{align*}
\pi_{*}= \pi_{*}\circ ( \phi \times \cdots \times \phi \times \id )^*: \Omega(\Delta_n \times X) \to \Omega(X),
\end{align*}
where $\pi$ denotes the projection $\Delta_n \times X \to X$.
This equation can be deduced inductively from the behavior of one-dimensional integrals under substitution.
\end{proof}

\begin{remark}
Later on, we will use the reparametrization invariance of differential forms in the image of Chen's map $\Chen$
under piecewise linear reparametrizations. In order to handle these kind of reparametrizations correctly, one should introduce the space of 
piecewise smooth paths $\tilde{\Path} M$
of $M$. By definition, a smooth map from $X$
 to $\tilde{\Path} M$ is a continuous maps $F^{\ev}: I \times X \to M$ together with a finite partition $([a_i,b_i])$ of 
$I$ such that
the restriction of $F^{\ev}$ to $X \times ]a_i,b_i[$ is smooth. The set $\smooth(X,\tilde{\Path} M)$ is acted upon by
the monoid of piecewise smooth reparametrizations of $I$ (still monotone and endpoints-preserving).
Observe that the definition of Chen's map $\Chen$ extends to smooth maps $X \to \tilde{\Path} M$ and still yields smooth
differential forms on $X$. Furthermore the reparametrization invariance continues to hold in the piecewise smooth setting.
\end{remark}



\subsection{Igusa's map} \label{subsection:Igusa_map}

In this paragraph we construct a map

\[\mathsf{S} :\Omega(\Path M) \rightarrow \s C(M)\]
from the differential forms on $\Path M$ to the suspension of the singular cochains on $M$.
We will need a sequence of maps from the cubes to the simplices, which relate their
cellular structures. These maps were originally considered by Adams (\cite{Adams}) and -- later on -- by
Chen  (\cite{C}). We will use a construction due
to Igusa (\cite{I}).

\begin{remark}
As mentioned before, we use the following definition of the $k$-simplex $\Delta_k$:
\begin{align*}
\Delta_k:=\{(t_1,\dots,t_k)\in \mathbb{R}^{k}: 1 \ge t_1 \ge t_2 \ge \cdots \ge t_k \ge 0\} \subset \mathbb{R}^{k}.
\end{align*}
In this convention, the face- and degeneracy-maps that equip $\{\Delta_k\}$ with the structure of a cosimplicial set are given by
\begin{eqnarray*}
\partial_i: \Delta_k \rightarrow \Delta_{k+1}, \quad (t_1,\dots, t_k) \mapsto
\begin{cases}
(1,t_1,\dots, t_k) & \textrm{for } i=0,\\
(t_1,\dots, t_{i-1}, t_i, t_i, t_{i+1}, \dots ,t_k) & \textrm{for } 0<i<k+1,\\
(t_1,\dots, t_k, 0) & \textrm{for } i=k+1.
\end{cases}
\quad \textrm{ and}\\
\epsilon_i: \Delta_k \rightarrow \Delta_{k-1}, \quad (t_1,\dots, t_k) \mapsto (t_1,\dots ,t_{i-1}, \widehat{t_i}, t_{i+1}, \dots, t_k),
\quad \textrm{respectively.}
\end{eqnarray*}
The $i$-th vertex $v_i$ of $\Delta_k$ is the point
\begin{align*}
(\underbrace{1, \dots, 1}_{\textrm{$i$ times}} ,\underbrace{0 ,\dots, 0}_{\textrm{$k-i$ times}}) \in \Delta_k.
\end{align*}

The simplicial set of (smooth) singular chains $\textrm{Sing}(M)$ of $M$ is given by
\begin{align*}
\textrm{Sing}_k(M) := \smooth(\Delta_k, M).
\end{align*} 
The simplicial structure maps $\{d_i\}$ and $\{s_i\}$ are defined via $d_i:=\partial_i^{*}$ and $s_i:=\epsilon_i^{*}$.

Observe that the maps $P_i$ and $Q_i$, which send an element $\textrm{Sing}(M)$ to its back- and front-face,  respectively, are equal to the pull backs of
\begin{eqnarray*}
U_i: \Delta_i &\to& \Delta_k, \quad (t_1,\dots, t_i) \mapsto (1,\dots,1,t_1,\dots,t_i) \quad \textrm{and}\\
V_i: \Delta_i & \to & \Delta_k, \quad (t_1,\dots,t_i)\mapsto (t_1,\dots,t_i,0,\dots,0), \quad \textrm{respectively.}
\end{eqnarray*}
\end{remark}

\begin{definition}
Let $M$ be a manifold.
The $dg$-algebra of (smooth) singular cochains $(C(M),\delta,\cup)$ consists of:
\begin{itemize}
\item The graded vector space $C(M)$ of linear functionals on the vector space generated by $\textrm{Sing}(M)$ over $\mathbb{R}$.
\item The differential $\delta$ defined by
\begin{align*}
(\delta \phi)(\sigma) := \sum_{i=0}^{k}(-1)^{i} (d_i^{*} \phi)(\sigma) := \sum_{i=0}^{k}(-1)^{i}\phi(\partial_i^{*}\sigma).
\end{align*}
\item The cup product $\cup$ defined by
\begin{align*}
(\phi \cup \psi)(\sigma) := \phi(V_i^{*} \sigma) \psi(U_j^{*} \sigma).
\end{align*}
\end{itemize}
\end{definition}

\begin{definition}
For each $k\geq 0$, the map
\[\Theta_{(k)}: I^{k-1} \rightarrow \Path \Delta_k, \]
 is defined to be the composition
 \begin{align*}
\xymatrix{
I^{k-1} \ar[r]^{\lambda_{(k)}}& \Path I^k\ar[r]^{\Path \pi_k}& \Path \Delta_k.
}
\end{align*}

Here $\pi_k: I^k \rightarrow \Delta_k$ is given by $\pi_k(x_1,\dots,x_k):=(t_1,\dots,t_k)$, with components
\[t_i:= \mathrm{ max }\{x_i,\dots, x_k\}.\]
The map $\lambda_{(k)}:I^{k-1}\rightarrow \Path I^k$
is defined by sending a point $(x_1,\dots,x_{k-1})$ to the path which goes backwards through the following $k+1$ points:
\[0 \leftarrow  x_1 e_1\leftarrow  (x_1 e_1+ x_2 e_2) \leftarrow \dots \leftarrow  ( x_1e_1+\dots + x_{k-1}e_{k-1})\leftarrow  (x_1e_1+\dots + x_{k-1}e_{k-1}+e_k), \]
where $(e_1,\cdots, e_n)$ denotes the standard basis of $\mathbb{R}^{n}$.
In other words, for $j=0,\dots, k$ we set:
\[\lambda_{(k)}(x_1,\dots, x_{k-1})(\frac{k-j}{k})= x_1 e_1+ \dots +x_{j} e_{j},   \]
where $x_k=1$, and interpolate linearly.

By convention, $\Theta_{(0)}$ is the map from a point to a point.

We will denote the map adjoint to $\Theta_{(k)}$ by $\Theta_k: I^k \rightarrow \Delta_k$.
\end{definition}

\begin{remark}

\end{remark}
Let us mention some properties of the maps $\Theta_{(k)}$:
\begin{enumerate}
\item The maps $\Theta_{(k)}$ are piecewise linear but not smooth. However, one can smooth them via
reparametrizations. The reparametrization invariance of elements in the image of Chen's map then guarantees
that for all our purposes the maps $\Theta_{(k)}$ behave as if they were smooth. 
\item One can check that the degree of the map $\Theta_k: I^{k} \to \Delta_k$ is $(-1)^{k}$, i.e.
\begin{align*}
\int_{I^{k}}\Theta_k^{*} \alpha = (-1)^{k} \int_{\Delta_k}\alpha
\end{align*} 
for every differential form $\alpha \in \Omega(\Delta_k)$.
\item The image of $\Theta_k$ lies in the subset $\Path(\Delta_k, v_k,v_0)$ of paths starting at the last vertex
$v_k$ and ending at the zeroth vertex $v_0$.
\end{enumerate}
The key properties of the maps $\Theta_{(k)}$ are described in the following lemma from \cite{I}:
\

\begin{lemma}
\label{lemma:commutative_diagrams}
The sequence of maps
\begin{align*}
\Theta_{(k)}: I^{k-1} \to \Path \Delta_k,
\end{align*}
satisfies the following properties:
\begin{enumerate}
\item For $1\leq i \leq k-1$ the diagram
\begin{align*}
\xymatrix{
I^{k-2} \ar[d]_{\Theta_{(k-1)}} \ar[r]^{\partial_i^-} & I^{k-1} \ar[r]^(.35){\Theta_{(k)}}& \Path(\Delta_k,v_k,v_0)\\
\Path(\Delta_{k-1},v_{k-1},v_0) \ar[rr]^{\hat{\phi}_i} &&  \Path(\Delta_{k-1},v_{k-1},v_0) \ar[u]_{\Path \partial_i}
}
\end{align*}
commutes,
where the maps are as follows:
\begin{itemize}
\item $\partial_i^-: I^{k-2}\to I^{k-1}$ is $(x_1,\cdots,x_{k-2}) \mapsto (x_1,\cdots,x_{i-1},0,x_i,\cdots,x_{k-2})$,
\item $\hat{\phi}_i$ is the map induced by the following piecewise smooth reparametrization
\begin{align*}
\phi_i(t):=
\begin{cases}
\frac{kt}{k-1} \quad &\textrm{for } 0\le t \le \frac{k-i-1}{k},\\
\frac{k-i-1}{k-1} \quad &\textrm{for } \frac{k-i-1}{k} \le t \le \frac{k-i}{k},\\
\frac{kt -1}{k-1} \quad &\textrm{for } \frac{k-i}{k} \le t \le 1.
\end{cases}
\end{align*}
Hence, $\hat{\phi}_i\circ\Theta_{(k-1)}=\phi_i\bullet \Theta_{(k-1)}$.
\item $\Path \partial_i$ is the smooth map induced by the the $i$'th face map $\partial_i:\Delta_{k-1} \to \Delta_{k}.$
\end{itemize}

\item  
For $1\leq i \leq k-1$ the diagram

\begin{align*}
\xymatrix{
I^{k-2} \ar[d]_{\cong} \ar[r]^{\partial_i^+} & I^{k-1} \ar[r]^(.4){\Theta_{(k)}}& \Path(\Delta_k,v_k,v_0)\\
I^{i-1}\times I^{k-i-1} \ar[rr]^(0.37){\Theta_{(i)} \times \Theta_{(k-i)}} &&  \Path(\Delta_{i},v_{i},v_0)\times \Path(\Delta_{k-i},v_{k-i},v_0) \ar[u]_{\mu_i},
}
\end{align*}
commutes.
Here $\mu_i$ is the concatenation map
\begin{align*}
\mu_i(\alpha,\beta)(t):=
\begin{cases}
U_{k-i}\left(\beta(\frac{kt}{k-i})\right) & \textrm{for } 0\le t \le \frac{k-i}{k},\\
V_i \left(\alpha(\frac{k}{i}(t-\frac{k-i}{k})\right) &\textrm{for } \frac{k-i}{k} \le t \le 1,
\end{cases}
\end{align*}
and
$\partial_i^+: I^{k-2}\to I^{k-1}$ is given by
 \[(x_1,\cdots,x_{k-2}) \mapsto (x_1,\cdots,x_{i-1},1,x_i,\cdots,x_{k-2}).\]
\end{enumerate}
\end{lemma}

\begin{lemma}\label{lemma:factorization}
Let $\mu_i$ be the $i$-th concatenation map from Lemma \ref{lemma:commutative_diagrams}
and suppose $a_1,\dots,a_n$ are differential forms on $\Delta_k$.
Also, let $f: X \to \Path(\Delta_i,v_i,v_0)$ and $g: Y \to \Path(\Delta_{k-i},v_{k-i},v_0)$
be smooth maps.
Then the following factorization property holds:
\begin{eqnarray*}
&& \hspace{-0.5cm} \int_{X\times Y} (f\times g)^{*} (\mu_i)^{*}\Chen(\s a_1 \otimes \cdots \otimes \s a_n) \\
&&= \sum_{l=0}^{n}\left(\int_X f^{*}\Chen(\s V_i^{*}a_1 \otimes \cdots \otimes \s V_i^{*}a_l)\right) \times
\left( \int_Y g^{*}\Chen(\s U_{k-i}^{*}a_{l+1}\otimes \cdots \otimes \s U_{k-i}^{*}a_n)\right).
\end{eqnarray*}

Here we extended the map $\Chen$ to the augmented bar complex $\mathbb{R}\oplus \B(\s \Omega(M))$ by setting
\begin{align*}
\Chen(\mathsf{1}) := 1.
\end{align*}
\end{lemma}

\begin{proof}
The main observation is that the diagram 
\begin{align*}
\xymatrix{
\Delta_l \times \Delta_{n-l} \times X \times Y \ar[rr]^{\cong}_{\tau}  \ar[d]_{\gamma_l \times (\id_{X\times Y})} && (\Delta_l \times X) \times (\Delta_{n-l}\times Y) \ar[d]^{(f^{\ev}) \times (g^{\ev})}\\
\Delta_n \times X\times Y \ar[rd]_{(\mu_i\circ (f\times g))^{\ev} \quad}&& (\Delta_i)^{l}\times (\Delta_{k-i})^{n-l} \ar[ld]^{\,\,\,\,\,\,\,\,(V_i)^{l}\times (U_{k-i})^{n-l}}\\
 &  (\Delta_k)^{n} &
}
\end{align*}
is commutative, where $\gamma_l$ is
\begin{eqnarray*}
\gamma_l: \Delta_l \times \Delta_{n-l} & \to & \Delta_n\\
(\vec{u},\vec{v}) &\mapsto & (\frac{i}{k}\vec{u} + \frac{k-i}{k}\vec{1}, \frac{k-i}{k}\vec{v}).
\end{eqnarray*}
It is easy to check that $\gamma_l$ is injective in the interior of $\Delta_l\times \Delta_{n-l}$ and maps onto
\begin{align*}
\Delta_n(l):=\{1\ge t_1 \cdots \ge t_l \ge \frac{k-i}{k} \ge t_{l+1} \ge \cdots \ge t_n \ge 0\}.
\end{align*}
Observe that the interiors of $\Delta_n(l)$ are disjoint for different values of $l$ and that
the union over all $l$ differs from $\Delta_n$ by a set of measure zero only.
Moreover, $\gamma_l$ is orientation preserving.

Now, we compute
\begin{eqnarray*}
&&\int_{X\times Y}(\pi_{X\times Y})_{*}((\mu_i\circ (f\times g)^{\ev})_{(n)})^{*}(p_1^{*}a_1\wedge \cdots \wedge p_n^{*} a_n)\\
&& = \int_{\Delta_n \times X\times Y} ((\mu_i\circ (f\times g)^{\ev})_{(n)})^{*}(p_1^{*}a_1\wedge \cdots \wedge p_n^{*} a_n)\\
&& = \sum_{l=0}^{n} \int_{\Delta_n(l)\times X\times Y}((\mu_i\circ (f\times g)^{\ev})_{(n)})^{*}(p_1^{*} a_1\wedge \cdots \wedge p_n^{*} a_n)\\
&& = \sum_{l=0}^{n} \int_{\Delta_l\times \Delta_{n-l} \times X\times Y}((\gamma_l \times\id_{X\times Y}))^{*}((\mu_i\circ (f\times g)^{\ev})_{(n)})^{*}(p_1^{*}a_1\wedge \cdots \wedge p_n^{*} a_n)\\
&& = \sum_{l=0}^{n}(-1)^{(n-l)\dim X} \int_{(\Delta_l \times X)\times (\Delta_{n-l} \times Y)}
((f^{\ev})_{(l)}\times (g^{\ev})_{(n-l)})^{*}\Big(\\
&& \hspace{6.5cm} ((V_i)^{l}\times (U_{k-i})^{n-l})^{*}
(p_1^{*}a_1\wedge \cdots \wedge p_n^{*} a_n)\Big)\\
&&  = \sum_{l=0}^{n}(-1)^{(n-l)\dim X} \left(\int_{\Delta_l \times X}(f^{\ev}_{(l)})^{*}(p_1^{*}V_i^*a_1\wedge \cdots \wedge p_l^*V_i^{*}a_l) \right) \times\\
&& \hspace{3.5cm}
\left(\int_{\Delta_{n-l}\times Y}(g^{\ev}_{(n-l)})^{*}(p_1^{*}U_{k-i}^*a_{l+1}\wedge \cdots \wedge p_{n-l}^*U_{k-i}^{*}a_n) \right)\\
&& = \sum_{l=0}^{n} (-1)^{(n-l)\dim X} \left(\int_{X}(\pi_{X})_*(f^{\ev}_{(l)})^{*}(p_1^{*}V_i^{*}a_1\wedge
\cdots \wedge p_l^{*}V_i^{*}a_l)\right) \times\\
&& \hspace{3.5cm} \left(\int_{Y}(\pi_{Y})_*(g^{\ev}_{(n-l)})^{*}(p_1^{*}U_{k-i}^{*}a_{l+1}\wedge
\cdots \wedge p_{n-l}^{*}U^*_{k-i}a_n)\right).\\
\end{eqnarray*}
Taking the additional signs in the definition of Chen's map $\Chen$ into account yields
the claimed factorization identity.
\end{proof}

\begin{definition}
The map $\mathsf{S}: \Omega(\Path M) \to \s C(M)$ is the composition of
\begin{eqnarray*}
\Omega(\Path M) &\to& C(M),\\
\alpha &\mapsto& \left( \sigma \mapsto \int_{I ^{k-1}}(\Theta_{(k)})^{*}\Path \sigma^{*}\alpha\right)
\end{eqnarray*}
and $\s: C(M) \to \s C(M)$.
\end{definition}

\begin{remark}
Let us evaluate the image of $\Chen(\s a)$ under $\mathsf{S}$ on a $k$-simplex
$\sigma: \Delta_k \to M$:
\begin{eqnarray*}
\mathsf{S}(\Chen(\s a))(\sigma)  &=& \int_{I^{k-1}}(\Theta_{(k)})^{*}\Path \sigma^{*}\Chen(\s a)
=\int_{I^{k-1}}(\Theta_{(k)})^*\Chen(\s \sigma^*a) \\
&=& \int_{I^{k}} (\Theta_{k})^{*}(\sigma^* a)
= (-1)^{k} \int_{\Delta_k}\sigma^{*}a. 
\end{eqnarray*}
Observe that this is valid only for $k>0$, since $\Chen(\s f)=0$ for $f$ a smooth function.
\end{remark}

\subsection{An $\mathsf{A}_{\infty}$ version of de Rham's theorem}\label{subsection:the_quasi_isomorphism}

The aim of this paragraph is to show that the composition
\begin{align*}
\xymatrix{
\B(\s \Omega(M)) \ar[rr]^{\Chen} &&\Omega(\Path M) \ar[rr]^{\mathsf{S}}&& \s C(M)
}
\end{align*}
yields an $\mathsf{A}_{\infty}$ quasi-isomorphism between the differential graded algebras $(\Omega(M),-d,\wedge)$
and $(C(M),\delta,\cup)$, see also \cite{Gugenheim}.

\begin{proposition}\label{proposition:Igusa}
Let $a_1,\dots, a_n$ be differential forms on $M$. Then the following equation holds:
\begin{eqnarray*}
\mathsf{S}(d\Chen(\s a_1\otimes \cdots \otimes \s a_n)) &=& {b}'_1(\mathsf{S}(\Chen(\s a_1 \otimes \cdots \otimes \s a_n))) +
\\ && \sum_{l=1}^{n-1}b'_2\left(\mathsf{S}(\Chen(\s a_1\otimes \s a_l))\otimes \mathsf{S}(\Chen(\s a_{l+1}\otimes \cdots \otimes \s a_n))\right). 
\end{eqnarray*}

Here $b'_1$ and $b'_2$ are the maps corresponding -- at the level of the suspension -- to the differential and the multiplication of the $dg$-algebra $C(M)$ of simplicial cochains.
\end{proposition}

\begin{proof}
Let $\alpha$ be an arbitrary differential form on the path space $\Path M$ and $\sigma: \Delta_k \to M$ a simplex. We want to compute
\begin{align*}
\int_{I^{k-1}}d(\Theta_{(k)})^{*}(\Path \sigma)^{*} \alpha = \int_{\partial I^{k-1}}\iota^{*}(\Theta_{(k)})^{*}(\Path \sigma)^{*}\alpha.
\end{align*}
Let $\partial_i^{\pm}$ be the standard embedings of $I^{k-2}$ into $I^{k-1}$ as top and bottom faces.
Then the expression above is equal to:
\begin{align*}
\sum_{i=1}^{k-1} (-1)^{i} \left(\int_{I^{k-2}}(\partial_{i}^{-})^{*}(\Theta_{(k)})^{*}(\Path \sigma)^{*}\alpha\right) - \sum_{i=1}^{k-1} (-1)^{i} \left(\int_{I^{k-2}}(\partial_{i}^{+})^{*}(\Theta_{(k)})^{*}(\Path \sigma)^{*}\alpha\right).
\end{align*}
We now take $\alpha = \Chen(\s a_1\otimes \cdots \otimes \s a_n)$.
Using the commutativity the first diagram  in Lemma \ref{lemma:commutative_diagrams},  reparametrization invariance of differential forms in the image of Chen's map $\Chen$ and
the naturality of Chen's map, one concludes
\begin{align*}
\int_{I^{k-2}}(\partial_i^-)^*(\Theta_{(k)})^{*}(\Path \sigma)^{*}\alpha = \int_{I^{k-2}} (\Theta_{(k-1)})^{*}(\Path \partial_i^{*}\sigma)^{*}\alpha=\mathsf{S}(\alpha)(\partial_i^{*}\sigma) .
\end{align*}
Consequently, the strata of $I^{k-1}$ given by the faces $\partial_i ^{-}I^{k-1}$ yield
\begin{align*}
\sum_{i=1}^{k-1}(-1)^i(\mathsf{S}(\Chen(\s a_1 \otimes \cdots \otimes \s a_n)))(\partial_i^{*}\sigma).
\end{align*}

On the other hand, commutativity of the second diagram  in Lemma \ref{lemma:commutative_diagrams} together with Proposition \ref{lemma:factorization} imply
\begin{eqnarray*}
\int_{I^{k-2}}(\partial^{+}_i)^{*}(\Theta_{(k)})^{*}(\Path \sigma)^{*}\alpha = 
\sum_{l=0}^{n}\mathsf{S}(\Chen(\s a_1\otimes \cdots \otimes \s a_l))(V_i^{*}\sigma) \mathsf{S}(\Chen(\s a_{l+1}\otimes \cdots \otimes \s a_n))(U_{k-i}^{*}\sigma).
\end{eqnarray*}
Summing over all the values of $i$ one obtains:
\begin{eqnarray*}
&& \hspace{-1cm} \sum_{i=1}^{k-1}(-1)^{i-1}\left(\int_{I^{k-2}}(\partial^{+}_i)^{*}(\Theta_{(k)})^{*}(\Path \sigma)^{*}\alpha\right)
=\\
&=&\sum_{l=1}^{n-1}b'_2(\mathsf{S}(\Chen(\s a_1\otimes \cdots \otimes \s a_l)) \otimes \mathsf{S}(\Chen(\s a_{l+1}\otimes \cdots \otimes \s a_n)))(\sigma)\\
&&+ \mathsf{S}(\Chen(\s a_1\otimes \cdots \otimes \s a_n))(U_{k-1}^{*}\sigma)+(-1)^k \mathsf{S}(\Chen(\s a_1\otimes \cdots \otimes \s a_n))(V_{k-1}^{*}\sigma).\\
&=&\sum_{l=1}^{n-1}b'_2(\mathsf{S}(\Chen(\s a_1\otimes \cdots \otimes \s a_l)) \otimes \mathsf{S}(\Chen(\s a_{l+1}\otimes \cdots \otimes \s a_n)))(\sigma)\\
&&+ \mathsf{S}(\Chen(\s a_1\otimes \cdots \otimes \s a_l))(\partial_0^{*} \sigma)+(-1)^k \mathsf{S}(\Chen(\s a_1\otimes \cdots \otimes \s a_l))(\partial_k^{*}\sigma).\\
\end{eqnarray*}
\end{proof}

\begin{definition}
Given a smooth manifold $M$ and $n\geq 1$, we define the map
\begin{align*}
\psi_n: \s \Omega(M)^{\otimes n} \to \s C(M)
\end{align*}
as follows:
\begin{enumerate}
\item For $n=1$, we set
\begin{align*}
\left(\psi_1(\s a)\right) (\sigma:\Delta_k \to M):= (-1)^{k} \left(\int_{\Delta^{k}}\sigma^{*}a \right).
\end{align*}
\item For $n>1$, we set
\begin{align*}
\psi_n(\s a_1 \otimes \cdots \otimes \s a_n) := (\mathsf{S}\circ \Chen)(\s a_1\otimes \cdots \otimes \s a_n).
\end{align*} 
\end{enumerate}
\end{definition}

\begin{remark}
Observe that $\psi_1(\s a)$ coincides with $(\mathsf{S}\circ \Chen)(\s a)$, except for the case
when $a$ is of degree $0$, i.e. a function. In that case, $(\mathsf{S} \circ \Chen)(\s a) = 0$, while
\begin{align*}
\left(\psi_1(\s a)\right)(\sigma: \{*\} \to M) := a(\sigma(0)).
\end{align*}
\end{remark}

We now come to the main result of this section, originally established  by Gugenheim in \cite{Gugenheim}.

\begin{theorem}[Gugenheim]\label{theorem:A_infty_quasi-isomorphism}
The sequence of maps:
\begin{eqnarray*}
\psi_n: (\s \Omega(M))^{\otimes n} \to \s C(M)
\end{eqnarray*}
gives an $\mathsf{A}_\infty$ morphism from $(\Omega(M),-d,\wedge)$ to $(C(M),\delta,\cup)$. Moreover, this morphism is a quasi-isomorphism and the construction is natural with respect to pull backs along smooth maps.
\end{theorem}

\begin{proof}
Let $\s a_1 \otimes \cdots \otimes \s a_n$ be an element of the bar complex
of the suspension of $(\Omega(M),-d,\wedge)$. By Theorem \ref{theorem:Chen}, the following equation holds:
\begin{eqnarray*}
d (\Chen(\s a_1\otimes \cdots \otimes  \s a_n)) &=& \Chen(\overline{D}(\s a_1\otimes \cdots \otimes \s a_n)) + \\
&& \hspace{-1cm} ev_1^{*}(a_1) \wedge \Chen(\s a_2 \otimes \cdots \otimes \s a_n) - (-1)^{[a_1]+\cdots + [a_{n-1}]} \Chen(\s a_1\otimes \cdots \otimes \s a_{n-1}) \wedge \ev_0^{*}(a_n).
\end{eqnarray*}
On the other hand, Proposition \ref{proposition:Igusa} asserts that
\begin{eqnarray*}
\mathsf{S}(d\Chen(\s a_1\otimes \cdots \otimes \s a_n) )&=& b'_1((\mathsf{S}\circ\Chen)(\s a_1 \otimes \cdots \otimes \s a_n)) +
\\ && \sum_{l=1}^{n-1}b'_2\left((\mathsf{S}\circ \Chen)(\s a_1\otimes\cdots \otimes \s a_l)\otimes (\mathsf{S}\circ\Chen)(\s a_{l+1}\otimes \cdots \otimes \s a_n)\right). 
\end{eqnarray*}
Together, these two equations yield
\begin{eqnarray*}
(\mathsf{S} \circ \Chen)(\overline{D}(\s a_1\otimes \cdots \otimes \s a_n)) &= &b'_1((\mathsf{S}\circ \Chen)(\s a_1\otimes \cdots \otimes \s a_n))\\
&& + \sum_{l=1}^{n-1}b'_2\left((\mathsf{S}\circ \Chen)(\s a_1\otimes \cdots \otimes \s a_l)\otimes (\mathsf{S}\circ \Chen)(\s a_{l+1}\otimes \cdots \otimes \s a_{n})\right)\\
&& - \mathsf{S}\left(\ev_1^{*}(a_1) \wedge \Chen(\s a_2 \otimes \cdots \otimes \s a_n) \right)\\
&& + (-1)^{[a_1]+\cdots + [a_{n-1}]}\mathsf{S}\left(\Chen(\s a_1\otimes \cdots \otimes \s a_{n-1}) \wedge \ev_0^{*}(a_0)\right).
\end{eqnarray*}
Direct computations lead to
\begin{align*}
\mathsf{S}\left(\ev_1^{*}(a_1)\wedge \Chen(\s a_2 \otimes \cdots \otimes \s a_n)\right) = -b'_2\left(\psi_1(\s a_1) \otimes (\mathsf{S}\circ \Chen)(\s a_2\otimes \cdots \otimes \s a_n) \right)
\end{align*}
for $|a_1|=0$
and
\begin{eqnarray*}
&&\mathsf{S}(\left(\Chen(\s a_1\otimes \cdots \otimes \s a_{n-1}) \wedge \ev_0^{*}(a_n)\right) =\\
&& \hspace{5cm} (-1)^{[a_{1}]+\cdots+[a_{n-1}]}b'_2\left((\mathsf{S}\circ \Chen)(\s a_2\otimes \cdots \otimes \s a_{n-1}) \otimes \psi_1(\s a_n) \right)
\end{eqnarray*}
for $|a_n|=0$, respectively. For $|a_1|>0$ (respectively, $|a_n|>0$), the first (second) expression vanishes.
Thus, we obtain:
\begin{eqnarray*}
&&\sum_{i=1}^{n}(-1)^{[a_1]+\cdots + [a_{i-1}]}\psi_{n}(\s a_1\otimes \cdots \otimes \s a_{i-1} \otimes \s (-da_i) \otimes \s a_{i+1} \otimes \cdots \otimes \s a_n) +\\
&&\sum_{i=1}^{n-1}(-1)^{[a_1]+ \cdots + [a_i]}\psi_{n-1}(\s a_1\otimes \cdots \otimes \s a_{i-1} \otimes \s (a_i \wedge a_{i+1}) \otimes \s a_{i+2} \otimes \cdots \otimes \s a_n) 
\\
&&  \hspace{3.5cm} = (\mathsf{S} \circ \Chen)(\overline{D}(\s a_1\otimes \cdots \otimes \s a_n)) = 
\\
&& b'_1(\psi_n(\s a_1\otimes \cdots \otimes \s a_n)) + \sum_{l=1}^{n-1}b'_2\left(\psi_l(\s a_1\otimes \cdots \otimes \s a_l)\otimes \psi_{n-l}(\s a_{l+1}\otimes \cdots \otimes \s a_{n})\right),\\
\end{eqnarray*}
which is precisely the structure equation for an $\mathsf{A}_{\infty}$ morphism.

We remark that, strictly speaking, the above argument is not valid for the case $n=2$ with both $a_1$ and $a_2$
of degree $0$, i.e. smooth functions. However, in this case the defining relation for $\psi$ being an $\mathsf{A}_{\infty}$ morphism
is equivalent to the fact that $(a_1 a_2)(x) = a_1(x) a_2(x)$ for every point $x\in M$.

That fact that the linear component
\begin{eqnarray*}
\psi_1: (\s \Omega(M))^{k-1} &\to& (\s C(M))^{k-1},\\
\alpha & \mapsto & (\sigma \mapsto (-1)^{k}\left(\int_{\Delta^k}\sigma^{*}\alpha \right))
\end{eqnarray*}
of the $\mathsf{A}_{\infty}$ morphism $\psi$ induces an isomorphism in cohomology is the content of the usual de Rham theorem.
Naturality of the construction follows from the naturality of $\mathsf{S}$ and $\Chen$.
\end{proof}

\begin{proposition}\label{lemma:some_properties}
The $\mathsf{A}_{\infty}$ quasi-isomorphism $\psi$ has the following properties:
\begin{enumerate}
\item $\psi_1(\s f) = \s f$ for every function $f\in  \Omega^{0}(M)$.
\item For $n>1$, $\psi_n(\s a_1\otimes \cdots \otimes \s a_n)$ vanishes, whenever one of the elements $a_i$ is a function.
\item The image of $\psi_n$ lies in the subspace $\hat{C}(M)$ of normalized simplicial cochains, i.e. the
space of cochains that vanish on degenerate simplices.
\end{enumerate}
\end{proposition}

\begin{proof}
The first two claims are direct consequences of the definition of $\psi$ and the properties of Chen's map $\Chen$.
Let us now prove the last claim. We need to show that 
\[ \psi_n(\s a_1 \otimes \dots \otimes \s a_n)(\epsilon_l^{*}\sigma)=0,\]
where $\epsilon_l: \Delta_k \rightarrow \Delta_{k-1}$ is the $l$-th degeneracy map.
Since
\begin{eqnarray*}
\psi_n(\s a_1 \otimes \dots \otimes \s a_n)(\epsilon_l^{*}\sigma)&=&\pm \int_{\Delta_n \times I^{k-1}} (\Theta_k)_{(n)}^*(\epsilon^{\times n}_l)^* (p^*_1\sigma^*(a_1)\wedge \dots \wedge p^*_n\sigma^*(a_n)),
\end{eqnarray*}
it is sufficient to show that the differential of
\begin{align*}
\xymatrix{
\Delta_n \times I^{k-1}\ar[rr]^(0.55){(\Theta_{k})_{(n)}}&&(\Delta_k)^{\times n}\ar[rr]^{{\epsilon_l}^{\times n}}&& (\Delta_{k-1})^{\times n}\\
}
\end{align*}
is singular almost everywhere.

Fix a point $(t_1,\dots,t_n,x_1,\dots,x_{k-1}) \in \Delta_n\times I^{k-1}$. We assume without loss of generality that each of the variables $t_m$ lies in an open interval of the form
\begin{align*}
\left(\frac{k-i_m}{k},\frac{k-i_m+1}{k}\right)
\end{align*}
where $i_m$ is an integer satisfying $1\le i_m \le k$.

Evaluating $(\Theta_k)_{(n)}$ on $(t_1,\dots,t_n,x_1,\dots,x_{k-1})$ yields an element of $(\Delta_k)^{n}$ whose projection onto the $m$-th copy is given by
\begin{eqnarray*}
\big(\max\{x_1,\dots,x_{i_m-1}, y_m \}, \dots, \max\{x_{i_m-1}, y_m \},
 y_m , \vec{0} \big),
\end{eqnarray*}
with $y_m:=k(t_m - \frac{k-i_m}{k})x_{i_m}$. We can assume without loss of generality that all the variables $(x_1,\dots,x_n,y_1,\dots,y_n)$ are pairwise different.

Now, consider the effect of applying the degeneracy-map $\epsilon_l: \Delta_k \to \Delta_{k-1}$ to all the components
of $(\Theta_k)_{(n)}(t_1,\dots,t_n,x_1,\dots,x_{k-1})$. The remaining expressions which (potentially) depend on $x_l$ and $(t_m)_{i_m=l}$ also depend on $x_{l-1}$, i.e.
we get either expressions of the form
\begin{eqnarray*}
\max\{x_j,\dots,x_{l-1},x_l,\dots\} \quad \textrm{or} \quad 
\max\{x_j,\dots, x_{l-1},k(t_m-\frac{k-l}{k})x_l,\vec{0}\}
\end{eqnarray*}
with $j\le l-1$,
or expressions independent of $x_l$ and $(t_m)_{i_m=l}$.
Hence, for 
\begin{align*}
((\epsilon_l)^{n} \circ (\Theta_k)_{(n)})(t_1,\dots,t_n,x_1,\cdots,x_n)
\end{align*}
to depend on $(x_l,(t_m)_{i_m=l})$, we need
\begin{align*}
x_{l-1} < k(t_m - \frac{k-l}{k}) x_l  \quad \textrm{for all } m \textrm{ with } i_m = l, \quad \textrm{and} \qquad x_{l-1} < x_l.
\end{align*}
Observe that this means that $x_{l-1}$ can only appear through expressions of the form
\begin{align*}
\max\{ \dots, k(t_m - \frac{k-l-1}{k})x_{l-1},\vec{0}\}
\end{align*}
where $i_m = l-1$. Therefore the restriction of the differential of $(\epsilon_l)^{n}\circ (\Theta_{k})_{(n)}$ at the point
$(t_1,\dots,t_n,x_{1},\dots,x_{k-1})$ to the subspace spanned by
\begin{align*}
\left\{\frac{\partial}{\partial x_{l-1}}, \left(\frac{\partial}{\partial t_m} \right)_{i_m=l-1} \right\}
\end{align*}
is not injective, and neither is the differential on the whole tangent space.
\end{proof}

\section{The integration $\mathsf{A}_{\infty}$ functor}\label{section:functor}

In this section the integration $\int[E]$ of a representation up to homotopy $E$
of a Lie algebroid $A$ is defined. The key concept is the holonomy map
\begin{align*}
\sigma \mapsto \Hol(\sigma,E),
\end{align*}
which assigns a linear map $\Hol(\sigma,E) \in \Hom^{1-k}(E_{\sigma(v_k)},E_{\sigma(v_0)})$
to every $k$-simplex $\sigma: T\Delta_k \to A$.

The main results are:
\begin{enumerate}
\item The assignment
\begin{eqnarray*}
\sigma \mapsto \Hol(\sigma,E)
\end{eqnarray*}
defines a unital representation up to homotopy of $\Pi_{\infty}(A)$ (Theorem \ref{theorem:integration_objects}).
\item The integration map
\begin{eqnarray*}
\int: \RRep(A) & \to & \URRep(\Pi_{\infty}(A)),\\
E & \mapsto & \Hol(-,E)
\end{eqnarray*}
can be extended naturally to an $\mathsf{A}_{\infty}$ functor of $dg$-categories (Theorem \ref{theorem:integration_everything}).
\end{enumerate}

Before proving these results, we need to establish one more property of the $\mathsf{A}_{\infty}$ quasi-isomorphism
\begin{align*}
\psi: (\Omega(M),-d,\wedge) \to (C(M),\delta,\cup).
\end{align*} 

\subsection{Gauge invariance}

Let $M$ be a smooth manifold and $V$ a finite dimensional graded vector space. We denote by ${\psi^V}$ the $\mathsf{A}_\infty$ morphism
\[ \psi^V:=\id_{\End V}\otimes  \psi: \End V\otimes (\Omega(M),-d,\wedge)\rightarrow \End V\otimes C(M). \]

In this paragraph, we will show that this $\mathsf{A}_\infty$ morphism is natural with respect to the gauge action. It turns out that this is a consequence of
general arguments regarding $\mathsf{A}_\infty$ morphisms between differential graded algebras. This fact will imply that the holonomies $\Hol(\sigma,E)$, associated to elements of $\Pi_{\infty}(A)$, are independent of the trivialization
of the graded vector bundle $\sigma^{*}E$.

\begin{remark}
Let $\mathsf{A}$ and $\mathsf{B}$ be  differential graded algebras with unit.
A Maurer-Cartan element $u\in \A^{1}$ gives rise to a differential graded algebra $\A_u$ with the same multiplication and twisted differential
\begin{align*}
d_u(a) := d a + u \wedge a - (-1)^{|a|} a \wedge u.
\end{align*}
Observe that this is a special instance of the twisting procedure of $\mathsf{A}_{\infty}$-algebras described in Appendix \ref{app2}.

Given an $\A_{\infty}$ morphism $\psi: \A \to \mathrm{B}$, we set
\[\psi(u):=\s^{-1} \left(\sum_{k\geq 1} \psi_k(\s u^{\otimes k})\right),\]
provided the sum converges. By Proposition \ref{MCgoestoMC}, $\psi(u)$ is a Maurer-Cartan element
of $\mathsf{B}$. 

As explained in the appendix \ref{app2}, there is an $\A_{\infty}$ morphism
\begin{align*}
\psi_u: \A_u \to \mathrm{B}_{\psi(u)},
\end{align*}
with structure maps $(\psi_u)_n$ given by
\[  (\psi_u)_n(\s a_1\otimes \dots \otimes \s a_n):=\sum_{l_0\ge 0, \cdots, l_n\ge 0} \psi_{n+l_0+\cdots+l_n}((\s u)^{\otimes l_0}\otimes \s a_1 \otimes \dots \otimes \s a_n \otimes (\s u)^{\otimes l_n}).\]
We assume convergence of this sum.
\end{remark}

\begin{definition}
Let $\A^{\times}$ be the group of invertible elements of $\A$ and $\textrm{MC}(\A)$
the set of Maurer-Cartan elements of $A$.
\end{definition}

\begin{remark}
Observe that there is a map
\begin{eqnarray*}
\A^{\times} &\to & \textrm{MC}(\A),\\
f &\mapsto & u_f:=f^{-1}df.
\end{eqnarray*}

Conjugation by $f\in \A^{\times}$
\begin{align*}
\phi_f(a):= f^{-1} a f
\end{align*}
yields an isomorphism of differential graded algebras
\begin{eqnarray*}
\phi_f: \A \to \A_{u_f}.
\end{eqnarray*}
\end{remark}

\begin{proposition}\label{propositiongauge}
 Let   $\psi:\A \to \B$
be an $\A_{\infty}$ morphism between two differential graded algebras with unit. 
Assume that $\psi$ satisfies the following properties:
\begin{enumerate}
\item  The map $\psi_1$  sends the unit to the unit.
\item For $n>1$, $\psi_n(\s a_1\otimes \cdots \otimes \s a_n)$ vanishes whenever one of the $a_i$ is invertible.
\end{enumerate}
Then, for any invertible element $f\in \A$, the following properties are satisfied:

\begin{enumerate}
\item The element $\psi_1(f)$  has an inverse given by $\psi_1(f^{-1})$.
\item $\psi(u_{f})$ equals $u_{\psi_1(f)} $.
\item The diagram
\begin{align*}
\xymatrix{
A \ar[r]^{\psi} \ar[d]_{\phi_f}^{\cong}& B \ar[d]^{\phi_{\psi_1(f)}}_{\cong} \\
A_{u_f} \ar[r]^{\psi_{u_f}} & B_{u_{\psi_1(f)}}
}
\end{align*}
is commutative.
\end{enumerate}

\end{proposition}

\begin{proof}
For the first claim we use the fact that $\psi$ is a morphism to compute:
\begin{eqnarray*}
0&=&b'_1 \psi_2(\s f,\s f^{-1})+b'_2(\psi_1(f)\otimes \psi_1(f^{-1}))-\psi_2(\s df\otimes \s f^{-1})+\psi_2(\s f\otimes \s d(f^{-1}))\\
&&-\psi_1(b_2(\s f \otimes \s f^{-1}))\\
&=&b'_2(\psi_1(f)\otimes \psi_1(f^{-1}))-\psi_1(b_2(\s f \otimes \s f^{-1}))=b'_2(\psi_1(f)\otimes \psi_1(f^{-1}))+\s 1.
\end{eqnarray*}

In order to establish the second claim, we first have to show that
\begin{align*}
u_{\psi_1(f)} = \psi_1(f)^{-1} d(\psi_1(f)) \qquad \textrm{and} \qquad \psi(u_f) = \sum_{k\ge 1}\psi_k((\s f^{-1}df)^{\otimes k})
\end{align*}
are equal. To this end, one observes that evaluating the structure equations for $\psi$ being an $\A_{\infty}$ morphism on the element
\begin{align*}
\s f \otimes (\s f^{-1}df)^{\otimes k}
\end{align*}
yields
\begin{align*}
b'_2(\psi_1(\s f)\otimes \psi_k((\s f^{-1}df)^{\otimes k})) =\psi_{k+1}(\s df \otimes (\s f^{-1}df)^{\otimes k})- \psi_k(\s df \otimes (\s f^{-1}df)^{\otimes k-1}) 
\end{align*}
and hence
\begin{align*}
\sum_{k\ge 1} b'_2(\psi_1(\s f)\otimes \psi_k((\s f^{-1}df)^{\otimes k})) =- \psi_1(\s df),
\end{align*}
which is equivalent to $u_{\psi_1(f)} = \psi(u_f)$.

In a similar manner one proves the following equations:
\begin{eqnarray*}
\sum_{l\ge 0} b'_2(\psi_1(\s f)\otimes \psi_{k+l+1}((sf^{-1}df)^{\otimes l}\otimes \s (f^{-1}a f) \otimes D)) &=& -\psi_{k+1}(\s (af) \otimes D),\\
\sum_{l\ge 0} \psi_{k'+k+l+2}(C\otimes \s (af) \otimes (s f^{-1}df)^{\otimes l} \otimes \s (f^{-1}bf) \otimes D)  & = & \psi_{k'+k+2}(C\otimes \s a \otimes \s (bf) \otimes D),\\
\sum_{l\ge 0} \psi_{k'+l+1}(C\otimes \s (af) \otimes (s f^{-1}df)^{\otimes l}) &=&(-1)^{[C]+[ a]} b'_2(\psi_{k'+1}(C\otimes \s a)\otimes \psi_1(\s f)),
\end{eqnarray*}
where $C\in (\s \A)^{k'}$ and $D\in (\s \A)^{k}$ are arbitrary homogeneous elements of 
the bar complex
\begin{align*}
\B(\s A)= \bigoplus_{k\ge 1}(\s A)^{\otimes k}.
\end{align*}

Together, these equations imply
\begin{eqnarray*}
&&  \sum_{l_0 \ge 0,\dots,l_n \ge 0} b'_2\left(   \psi_1(\s f)\otimes \psi_{l_0+\cdots + l_n + n}\left((\s f^{-1}df)^{\otimes l_0} \otimes  \s (f^{-1}a_1 f) \otimes (\s f^{-1}df)^{\otimes l_1} \otimes \cdots \right. \right.\\
&& \left. \left.\hspace{2 cm} \cdots \otimes \s (f^{-1}a_n f)\otimes (\s f^{-1}df)^{\otimes l_n}\right)\right) \\
&& =(-1)^{[a_1]+\dots+[a_n]+1}b'_2 \left(\psi_n(\s a_1\otimes \cdots \otimes \s a_n) \otimes \psi_1(\s f)\right),
\end{eqnarray*}
which is another way to express
\begin{align*}
\psi_{u_f}\circ \phi_f = \phi_{\psi_1(f)} \circ \psi.
\end{align*}
\end{proof}

\begin{definition}
The gauge action of a differential graded algebra $\A$ is the right action of the group of invertible elements $\A^{\times}$ on the set $\A$ given by
\[ a \bullet f:= f^{-1}af +f^{-1}df.  \]
\end{definition}

\begin{corollary}\label{corollary:gauge}
Let $M$ be a smooth manifold and $V$ be a finite dimensional graded vector space. We denote by ${\psi^V}$ the $\mathsf{A}_\infty$ morphism
\[ \psi^V:=\id_{\End V}\otimes \psi: \End V\otimes (\Omega(M),-d,\wedge)\rightarrow \End V\otimes C(M), \]
given by the tensor product of $\id_{\End V}$ and the $\mathsf{A}_{\infty}$ quasi-isomorphism $\psi$, see Appendix \ref{app2}.

Then:
\begin{enumerate}
\item For any $a\in \End V\otimes \Omega(M)$, the sequence 
\[\psi^{V}(a):= \sum_{k \geq 1} \s^{-1}(\psi^V_k((\s a)^{\otimes k}))\]
converges, in the sense that evaluating every term of the above sum on some fixed simplex yields a sum of endomorphisms which converges absolutely.
\item The map $a \mapsto \psi^E(a)$ is equivariant with respect to the gauge action, namely:
\[ \psi^E( a\bullet f)=\psi^E(a)\bullet f. \]
Here we use the fact that the invertible elements of $\End V\otimes (\Omega(M),-d,\wedge)$ can also be seen as invertible elements of $\End V\otimes C(M)$.
\end{enumerate}
\end{corollary}

\begin{proof}
The first claim follows from a simple bound that one obtains from the fact that the volume of the $n$-simplex is $\frac{1}{n!}$.
The second claim is a formal consequence of Proposition \ref{propositiongauge}, since Lemma \ref{lemma:some_properties} implies that $\psi^V$ satisfies 
the hypothesis there.
\end{proof}

\subsection{Integration of representations up to homotopy}\label{subsection:integration_objects}

In this paragraph we show that the $\mathsf{A}_\infty$ quasi-isomorphism $\psi$ from the de Rham algebra to the algebra of singular cochains, which was constructed in Section \ref{section:diff_forms_and_cochains},
can be used to integrate representations up to homotopy. 
In particular, we use $\psi$ to define the holonomy
\begin{align*}
\Hol(\sigma,E)
\end{align*}
of a simplex $\sigma \in \Pi_{\infty}(A)$ with respect to some representation up to homotopy $E$ of $A$.
Theorem \ref{theorem:integration_objects} asserts that, in fact, the assignment
\begin{align*}
\sigma \mapsto \Hol(\sigma,E)
\end{align*}
defines a representation up to homotopy of $\Pi_{\infty}(A)$.

Before that, we have to modify slightly the map $\psi$, by reversing the orientations
in order to be consistent with Igusa's conventions in \cite{I}, as well as the conventions
of \cite{AC2}\footnote{The need for this change of orientation comes from the fact that Igusa's maps
$\Theta_{(k)}$ are maps from $I^{k-1}$ to the space of paths in $\Delta_k$ from $v_k$ to $v_0$. This is consistent with the conventions  for the nerve of a category in \cite{AC2}, where by
a sequence of composable arrows $(g_1,\dots,g_k)$ we mean a sequence such that the target of $g_i$ is equal to the source of $g_{i-1}$.  The discrepancy of these conventions with the usual orientations of simplices and nerves forces us to introduce this change.}.

\begin{definition}

Given a manifold $M$, we denote by $\overline{C}(M)$ the $dg$-algebra $(C(M),  \overline{\delta},\overline{\cup})$ with:
\begin{itemize}
\item[i)] underlying vector space the space of smooth singular cochains $C(M)$,
\item[ii)] differential 
\[\overline{\delta} : C^k(M )\rightarrow C^{k+1}(M), \quad \overline{\delta}:=(-1)^{k}\delta \] 
and 
\item[iii)] product given by
\[\alpha \overline{\cup} \beta:=  (-1)^{|\alpha||\beta|} \alpha \cup \beta. \]
\end{itemize}
Let $\overline{i}$ be the natural isomorphism of differential graded algebras
\begin{eqnarray*}
\overline{i}: C(M)\rightarrow \overline{C}(M),\\
\alpha \mapsto (-1)^{\frac{|\alpha|(|\alpha|-1)}{2}}\alpha
\end{eqnarray*}
and $\tau$ the natural isomorphism of differential graded algebras
\begin{eqnarray*}
\tau: (\Omega(M), d,\wedge) \rightarrow (\Omega,-d,\wedge)  \\
\eta \mapsto (-1)^{|\eta |}\eta.
\end{eqnarray*}
We define

\[\overline{\psi}:=\overline{i}\circ \psi \circ \tau .\]
Clearly, $\overline{\psi}: (\Omega(M),d,\wedge) \rightarrow \overline{C}(M)$ is an $\mathsf{A}_\infty$ quasi-isomorphism. 
Moreover, its components $\overline{\psi}_n:\s\Omega(M)^{\otimes n}\rightarrow \s \overline{C}(M)$ are given by
\[\overline{\psi}_n(\s a_1 \otimes \dots \otimes \s a_n)=(-1)^{\frac{k(k-1)}{2}+k+n-1}\psi_{n}(\s a_1\otimes \dots \otimes \s a_n ),\]
where $k=[a_1]+\dots +[a_n]+1.$

\end{definition}

\begin{remark}
Let $M$ be a manifold and $V$ a finite dimensional graded vector space. 
We identify the vector space $ \End V\otimes  C(M)$ with the space of cochains on $M$ with values in
$\End V$ as follows. For an element $\phi \otimes \eta \in \End V\otimes  C^k(M)$ and a $k$-simplex $\sigma$
in $M$, we set
\[(\phi \otimes \eta)(\sigma):=\phi \eta (\sigma)\in \End V.\]  
A simple computation shows that if $\alpha \in \End^l V\otimes C^k(M)$ and $\beta \in \End^{l'} V \otimes C^{k'}(M)$ then: 
\[(\alpha \overline{\cup} \beta)(\sigma)=(-1)^{k(l'+k')}\alpha(V_k^{*}(\sigma)) \circ \beta (U_{k'}^{*}(\sigma)). \]
\end{remark}

\begin{proposition}\label{lemmaMC2}
Let $M$ be a smooth manifold and $V$ a finite dimensional graded vector space. 
\begin{itemize}
\item[a)]
There is a natural bijective correspondence
between
\begin{enumerate}
\item Maurer-Cartan elements in the $dg$-algebra $\End V\otimes \overline{C}(M)$.
\item Representations up to homotopy of the simplicial set $\Pi_\infty(M):=\textrm{Sing}(M)$, such that the graded vector space associated
to every point is $V$.
\end{enumerate}
The correspondence is given by
\[\alpha\mapsto \bold{1}+\alpha \in \End V\otimes C(M). \]
Here $\bold{1}$ is the one-cochain with values in $\End V$, which associates the identity to every $1$-simplex.

\item[b)]
Suppose that $E$ and $E'$ are representations up to homotopy of $\Pi_\infty(M)$ on trivial vector bundles $M\times V$ and $M\times V'$, respectively. Then there is a natural isomorphism
of vector spaces:
\begin{eqnarray*}
 \Hom(V,V')\otimes C(M)&\cong& \RHom(E,E'),\\
\end{eqnarray*}
where $\RHom(E,E')$ denotes the complex of morphisms between $E$ and $E'$.

Under this identification, the operator $\D:\RHom(E,E')\rightarrow \RHom(E,E')$ corresponds to the map
\begin{eqnarray*}
 \Hom(V,V')\otimes C(M)&\rightarrow &  \Hom(V,V')\otimes C(M) \\
 \eta &\mapsto& \overline{\delta} \eta +\alpha' \overline{\cup} \eta -(-1)^{|\eta|}\eta \overline{\cup} \alpha,
\end{eqnarray*}
where $\alpha$ and $\alpha'$ are the Maurer-Cartan elements corresponding to $E$ and $E'$ and the product
is taken in the algebra $\End(V\oplus V')\otimes \overline{C}(M)$.
Moreover, under the identification, composition corresponds to the product in the algebra $\End(V\oplus V' \oplus V'')\otimes \overline{C}(M)$.
\end{itemize}
\end{proposition}

\begin{proof}
We write $F=\bold{1}+\alpha$. The fact that $\alpha$ is a Maurer-Cartan element is equivalent to:

\[\overline{\delta}(F-\bold{1})+ (F-\bold{1}) \overline{\cup} (F-\bold{1})=0.\]
Since $-\bold{1}$ is a Maurer-Cartan element of $\End V\otimes \overline{C}(M) $, this equation becomes
\begin{equation}\label{eqwith1}
\overline{\delta}(F)+ (F) \overline{\cup} (F) - \bold{1}\overline{\cup} F - F \overline{\cup} \bold{1}=0. 
\end{equation}
We know that $F$ can be written as a sum of homogeneous components
\[F=F_0+F_1+F_2+\dots, \]
where $F_k\in \End^{1-k}(V)\otimes \overline{C}^k(M)$. Looking at the homogeneous components, equation \eqref{eqwith1} is equivalent to
the sequence of equations
\[ \overline{\delta}(F_{k-1})+ \sum_{i+j=k} F_j \overline{\cup} F_i  - \bold{1}\overline{\cup} F_{k-1} - F_{k-1} \overline{\cup} \bold{1}=0. \]
Using the fact that $\bold{1}\overline{\cup} F_{k-1}=(-1)^{k-1}(\id \otimes d^*_0)(F_{k-1}) $ and $F_{k-1} \overline{\cup} \bold{1}=-( \id \otimes d^*_k) (F_{k-1})$, we obtain that the equation above is equivalent to
\[ \sum_{i=1}^{k-1}(-1)^{k-1+j}(\id \otimes d^*_i)(F_{k-1})+ \sum_{i+j=k} F_j \overline{\cup} F_i =0. \]
This equations are precisely the defining relations of a representation up to homotopy. This concludes the proof of the first statement.

For the second statement, we take an element
 \[\eta=\eta_0+\eta_1+\eta_2+\dots\]
in $(\Hom(V, V')\otimes \overline{C}(M))^n$ and evaluate $\overline{\delta} \eta + \alpha' \overline{\cup} \eta -(-1)^n \eta \overline{\cup} \alpha$ on a $k$-simplex $\sigma$.
\begin{eqnarray*}
\left(\overline{\delta} \eta + \alpha' \overline{\cup} \eta -(-1)^n \eta \overline{\cup} \alpha \right)(\sigma)&=&\left(\overline{\delta} \eta_{k-1} +\sum_{i+j=k} \alpha'_j \overline{\cup} \eta_i -\sum_{i+j=k}(-1)^n \eta_j \overline{\cup} \alpha_i \right)(\sigma)\\
&=&\left( \sum_{i=1}^{k-1}(-1)^{k-1+i}(\id \otimes {d_i^*}) \eta_{k-1} +\sum_{i+j=k} F'_j \overline{\cup} \eta_i \right)(\sigma)\\
&&-\sum_{i+j=k}(-1)^n (\eta_j \overline{\cup} F_i )(\sigma)\\
&=&\sum_{j=1}^{k-1}(-1)^{j+n}\eta_{k-1} (d_j(\sigma))+\sum_{i+j=k}(-1)^{jn}F_{j}^{'}\cup
 \eta_{i}(\sigma)\\
 &&+ \sum_{i+j=k}(-1)^{n+j+1}\eta_j \cup
F_{i}(\sigma).
\end{eqnarray*}
The last expression is precisely the formula for the differential on the spaces of morphisms defined in Equation \eqref{differentialonmorphisms}.
The remaining statements follow from similar computations.
\end{proof}

\begin{definition}
Let $E$ be a representation up to homotopy of a Lie algebroid $A$.

The holonomy $\Hol(\sigma,E)$ of a simplex $\sigma: T\Delta_k \to A$ of $\Pi_{\infty}(A)$
is the linear map
\begin{align*}
\Hol(\sigma,E): E_{\sigma(v_k)} \to E_{\sigma(v_0)}
\end{align*} 
of degree $1-k$ defined as follows:
\begin{enumerate}
\item Pull back the representation up to homotopy along $\sigma$
and choose a trivialization
\begin{align*}
h: \sigma^{*}E \cong \Delta_k\times V
\end{align*}
 of the graded vector bundle $\sigma^{*}E \to \Delta_k$.
\item By Proposition \ref{lemmaMC1}, the representation up to homotopy on $\Delta_k \times V$ corresponds to a Maurer-Cartan element $\omega$ of $\End V \otimes \Omega(\Delta_k)$.
Apply the $\mathsf{A}_{\infty}$ quasi-isomorphism
\begin{align*}
\overline{\psi}^{V}: \End V \otimes \Omega(\Delta_k) \to \End V \otimes \overline{C}(\Delta_k)
\end{align*}
to $\omega$. This yields a cochain $\overline{\psi}^{V}(\omega)$ on $\Delta_k$ with
values in $\End V$.
\item Evaluating this cochain on the fundamental cycle $[\Delta_k]:=(\id: \Delta_k \to \Delta_k)$ gives an element of $\End V$, which we interpret as a linear map
from the fibre of $\Delta_k \times V$ over $v_k$ to the fibre over $v_0$. Finally, we set
\begin{align*}
\Hol(\sigma,E) := h_{v_0}^{-1} \circ <\overline{\psi}^{V}(\omega),[\Delta_k] > \circ h_{v_k}.
\end{align*}
\end{enumerate}
\end{definition}

\begin{lemma}\label{lemma:invariance}
The holonomy $\Hol(\sigma,E)$ of a simplex $\sigma \in \Pi_{\infty}(A)$
is well-defined.
\end{lemma}

\begin{proof}
We have to prove that the linear map $\Hol(\sigma,E)$ is independent of the chosen
trivialization $h: \sigma^{*}E \cong \Delta_k\times V$.
Changing the trivialization can be encoded in an automorphism of graded vector bundles
\begin{align*}
f: \Delta_k \times V \to \Delta_k \times V,
\end{align*}
which can also be seen as an invertible element of the $dg$-algebra $\End V \otimes \Omega(\Delta_k)$.

If one uses the trivialization $f^{-1}\circ h$, the Maurer-Cartan element $\omega$
changes to
\begin{align*}
\omega' = f^{-1} \omega f + f^{-1}df = \omega \bullet f.
\end{align*}
Using Corollary \ref{corollary:gauge}, we compute
\begin{eqnarray*}
<\bold{1} + \overline{\psi}^{V}(\omega \bullet f), [\Delta_k]> & = & <\bold{1} + \overline{\psi}^{V}(\omega) \bullet f, [\Delta_k]>\\
&=& <\bold{1} + f^{-1}\overline{\delta} f, [\Delta_k]> + <f^{-1}\overline{\cup} \overline{\psi}^{V}(\omega) \overline{\cup} f, [\Delta_k]>\\
&= & f^{-1}(v_0) \circ <\bold{1},[\Delta_k]> \circ f(v_k) + f^{-1}(v_0) \circ <\overline{\psi}^{V}(\omega),[\Delta_k]> \circ f(v_k) \\
& =& f^{-1}(v_0) \circ \left(<\bold{1} + \overline{\psi}^{V}(\omega),[\Delta_k]> \right) \circ f(v_k).
\end{eqnarray*}

This means that different choices of trivializations lead to linear maps, which are related by conjugation.
Consequently 
\begin{align*}
\Hol(\sigma,E):=  h_{v_o}^{-1} \circ <\overline{\psi}^{V}(\omega),[\Delta_k]> \circ h_{v_k}
\end{align*}
is independent of the trivialization.
\end{proof}

\begin{remark}\label{remark:parallel_transport}
\hspace{0cm}
\begin{enumerate}
\item The definition of the holonomy $\Hol(\sigma,-)$
of a simplex $\sigma: T\Delta_k \to A$ does not make use of all
the defining properties of a representation up to homotopy.
In fact, it only uses the fact that the structure operators of the pull back $\sigma^{*}E \cong \Delta_k\times V$
can be naturally assembled into a differential form $\omega$ with values in $\End V$.
This is possible because the linear operator
\begin{align*}
D: \Omega(A,E) \to \Omega(A,E)[1],
\end{align*}
which encodes the representation up to homotopy, satisfies the graded derivation rule
\begin{align*}
D(\omega \eta) = d(\omega) \eta + (-1)^{|\omega|}\omega D(\eta).
\end{align*}
Hence, our definition of holonomy can be extended to $\mathbb{Z}$-graded $A$-connections, i.e. linear operators $D$ of $\Omega(A,E)$ of degree $1$, which satisfy
the graded derivation rule, but do not necessarily square to zero.

\item Any connection $\nabla$  on a vector bundle $E\to M$ gives rise to a
$\mathbb{Z}$-graded $TM$-connection
\begin{align*}
d_{\nabla}: \Omega(A,E) \to \Omega(A,E)
\end{align*}
via the Chevalley-Eilenberg formula
\begin{eqnarray*}
d \eta(\alpha_1, \dots ,\alpha_{n+1})&=& \sum_{i<j}(-1)^{i+j} \eta([\alpha_i,\alpha_j],\cdots,\hat{\alpha}_i,\dots,\hat{\alpha}_j,\dots,\alpha_{k+1})\\
&&+\sum_i(-1)^{i+1}\nabla_{(\alpha_i)}\eta(\alpha_1,\dots,\hat{\alpha}_i,\dots,\alpha_{k+1}).
\end{eqnarray*}

We claim that the holonomy
\begin{align*}
\Hol(\gamma,E): E_{\gamma(1)} \to E_{\gamma(0)}
\end{align*}
associated to a path $\gamma$ is equal to the parallel transport of $\nabla$ along $t \mapsto \gamma(1-t)$.

To see this, we pick a trivialization $\gamma^{*}E \cong I \times V$ and write
\begin{align*}
\gamma^{*}\nabla  = d + a(t)dt.
\end{align*}
Evaluating the cochain
\begin{align*}
\overline{\psi}^{V}_{n}((\s a)^{n}) = (-1)^{n} \psi_n^{V}((\s a)^{n})
= (-1) ^{n} (\mathsf{S} \circ \Chen)((\s a)^{n})
\end{align*}
on the fundamental cycle $[I]=(\id: I \to I)$ yields
\begin{align*}
\int_{1\ge t_1 \ge \cdots \ge t_n \ge 0} a(1-t_1) \circ \cdots \circ a(1-t_n)dt_1 \cdots dt_n.
\end{align*}
Hence -- up to the identification $\sigma^{*}E \cong I \times V$ -- $\Hol(\gamma,E)$ is given by
\begin{align*}
H:=\id_{V} + \sum_{n\ge 1}\left(\int_{1\ge t_1 \ge \cdots \ge t_n \ge 0} a(1-t_1) \circ \cdots \circ a(1-t_n)dt_1 \cdots dt_n\right).
\end{align*}
Observe that the one-parameter family of endomorphisms
\begin{align*}
H_t:= \id_{V} + \sum_{n\ge 1}\left(\int_{t\ge t_1 \ge \cdots \ge t_n \ge 0} a(1-t_1) \circ \cdots \circ a(1-t_n)dt_1 \cdots dt_n\right)
\end{align*}
satisfies
\begin{align*}
H_0 = \id_{V} \quad \textrm{and} \quad \frac{d}{dt}H_t = a(1-t) \circ H_t,
\end{align*}
which is the differential equation defining the parallel transport of $\nabla$ along
$t \mapsto \gamma(1-t)$.
\end{enumerate}
\end{remark}

\begin{theorem}\label{theorem:integration_objects}
Let $A$ be a Lie algebroid and $E$ a representation up to homotopy of $A$.

The assignments
\begin{eqnarray*}
(\Pi_{\infty}(A))_0 \ni x &\mapsto& E_x,\\
\Pi_{\infty}(A) \ni \sigma &\mapsto &F_k(\sigma):= \Hol(\sigma,E)
\end{eqnarray*}
define a unital representation up to homotopy of the $\infty$-groupoid $\Pi_{\infty}(A)$
of $A$.
\end{theorem}

\begin{proof}
We can assume that $\sigma^{*}E$ is trivial, i.e. $\sigma^{*}E \cong \Delta_k \times V$.
Since $\omega$ is a Maurer-Cartan element of $\End V \otimes \Omega(\Delta_k)$,
so is $\overline{\psi}^{V}(\omega)$ for $\End V \otimes \overline{C}(\Delta_k)$.
By Proposition \ref{lemmaMC2}, such a Maurer-Cartan element corresponds uniquely to a representation
up to homotopy of $\Pi_{\infty}(\Delta_k) = \textrm{Sing}(\Delta_k)$.
The defining relations of such a representation up to homotopy yield the defining relations for $\Hol(-,E)$ being a representation up to homotopy, evaluated on the simplex $\sigma$. 

Unitality follows from part $3$ of Lemma \ref{lemma:some_properties}. 
\end{proof}

\begin{definition}
Let $A$ be a Lie algebroid. The integration map
\begin{align*}
\int: \RRep(A) \to \URRep(\Pi_{\infty}(A))
\end{align*}
assigns to any representation up to homotopy $E$ of $A$
the representation up to homotopy $\int[E]$ of $\Pi_{\infty}(A)$.
\end{definition}

\begin{remark}\label{remark:natural}
Every morphism of Lie algebroids $f: A \to B$ yields a $dg$-functor
\begin{align*}
f^{*}: \RRep(B) \to \RRep(A),
\end{align*}
as well as a morphism of simplicial sets
\begin{align*}
\Pi_{\infty}f: \Pi_{\infty}(A) \to \Pi_{\infty}(B),
\end{align*}
which, in turn, yields
\begin{align*}
(\Pi_{\infty}f)^{*}: \URRep(\Pi_{\infty}(B)) \to \URRep(\Pi_{\infty}(A)).
\end{align*}
The map $\int$ is natural with respect to pull back, i.e. the diagram
\begin{align*}
\xymatrix{
\URRep(\Pi_{\infty}(B)) \ar[rr]^{(\Pi_{\infty}f)^{*}} && \URRep(\Pi_{\infty}(A))\\
\RRep(B) \ar[rr]^{f^{*}} \ar[u]^{\int} & & \RRep(A) \ar[u]_{\int}
}
\end{align*}
commutes.
\end{remark}

\subsection{The $\mathsf{A}_{\infty}$ functor}

Next, we extend the map
\begin{eqnarray*}
\int: \RRep(A) &\to & \URRep(\Pi_{\infty}(A))\\
E & \mapsto & \int[E]:= \Hol(-,E),
\end{eqnarray*}
which was defined in the previous subsection, to an $\mathsf{A}_{\infty}$ functor
of $dg$-categories
\begin{align*}
\int: \RRep(A) \to \URRep(\Pi_{\infty}(A)).
\end{align*}
The structure maps
\begin{align*}
\int_n: \s \RHom(E_1,E_0) \otimes \cdots \otimes \s \RHom(E_n,E_{n-1}) \to \s \RHom(E_n,E_0)
\end{align*}
of this $\mathsf{A}_{\infty}$ functor
are given in terms of the holonomy $\Hol(\sigma, \phi_1,\dots, \phi_n)$
of a simplex $\sigma: T\Delta_k \to A$ with respect to
a chain of composable morphisms
\begin{align*}
\xymatrix{
E_0  & E_{1} \ar[l]_{\phi_1} & \cdots \ar[l]& E_{n-1} \ar[l] & E_n \ar[l]_{\phi_n}.
}
\end{align*}

\begin{definition}
Let $E_1,\dots, E_n$ be representations up to homotopy of a Lie algebroid $A$
and $\phi_i: E_i \to E_{i-1}$ morphisms in $\RRep(A)$.

The holonomy $\Hol(\sigma,\phi_1,\dots, \phi_n)$ of a simplex $\sigma: T\Delta_k \to A$
is the linear map
\begin{align*}
\Hol(\sigma,\phi_1,\dots, \phi_n): E_n|_{\sigma(v_k)} \to E_0|_{\sigma(v_0)} 
\end{align*}
of degree $[\phi_1]+\cdots+[\phi_n]-k+1$ defined as follows:
\begin{enumerate}
\item Pull back the representations up to homotopy $E_i$ and all the morphisms $\phi_i$ along $\sigma$
and choose a trivialization
\begin{align*}
h: \sigma^{*}(E_0 \oplus \cdots \oplus E_n) \cong \Delta_k \times (V_0\oplus \cdots \oplus V_n) =: \Delta_k \times V.
\end{align*}
\item By Proposition \ref{lemmaMC1}, the representations up to homotopy on $\Delta_k \times V_i$
corresponds to a Maurer-Cartan element $\omega^i$ of $\End V \otimes \Omega(\Delta_k)$,
while the morphisms correspond to elements $\eta^i \in \Hom(V_i,V_{i-1})\otimes \Omega(\Delta_k)$.
\item
Apply the $\mathsf{A}_{\infty}$ quasi-isomorphism
\begin{align*}
\overline{\psi}^{V}: \End V \otimes \Omega(\Delta_k) \to \End V \otimes \overline{C}(\Delta_k)
\end{align*}
to $\Omega:=\omega_0 + \cdots + \omega_n + \eta_1 + \cdots + \eta_n$.
This yields a cochain $\overline{\psi}^{V}(\Omega)$ on $\Delta_k$ with values in $\End V$.
\item
Evaluating this cochain on the fundamental cycle $[\Delta_k]$ gives an element in $\End V$.
We define $\Hol(\sigma, \phi_1,\dots,  \phi_n)$ to be the composition:
\begin{align*}
\xymatrix{
E_n|_{\sigma(v_k)} \ar[rr]^{h_{v_k}} && V \ar[rrr]^{<\overline{\psi}^{V}(\Omega),[\Delta_k]>}&&& V \ar[rr]^{h_{v_0}^{-1}}&& E_0|_{\sigma(v_0)}.
}
\end{align*}
\end{enumerate}
\end{definition}

\begin{lemma}\label{lemma:invariance2}
The holonomy $\Hol(\sigma, \phi_1,\dots, \phi_n)$ is well-defined.
\end{lemma}

\begin{proof}
The proof is literally the same as the one of Lemma \ref{lemma:invariance}.
\end{proof}

\begin{theorem}\label{theorem:integration_everything}
Let $A$ be a Lie algebroid.
The assignments
\begin{align*}
\RRep(A) \ni E  \mapsto \int[E]\in \URRep(\Pi_{\infty}(A))
\end{align*}
and
\begin{eqnarray*}
\int_n: \s \RHom(E_1,E_0) \otimes \cdots \otimes \s \RHom(E_n,E_{n-1}) &\to & \s \RHom(E_n,E_0),\\
\phi_1 \otimes \cdots \otimes \phi_n & \mapsto & \s \left(\Hol(-,\phi_1, \dots,  \phi_n)\right)
\end{eqnarray*}
define an $\mathsf{A}_{\infty}$ functor 
\begin{align*}
\int: \RRep(A) \to \URRep(\Pi_{\infty}(A))
\end{align*}
between the $dg$-category of representations
up to homotopy of $A$ and the $dg$-category of unital representations up to homotopy of $\Pi_{\infty}(A)$.
\end{theorem}

\begin{proof}
We need to prove the equations
\begin{equation}\label{equationtoprove}
 \beta'_1\circ \int_{n}+\sum_{i+j=n}\beta'_2\circ \left( \int_i \otimes \int_j\right)=\sum_{i+j+1=n}\int_{n}\circ \left(\id^{\otimes i} \otimes \beta_1\otimes \id^{\otimes j}\right)+ \sum_{i+j+2=n} \int_{n-1}\circ \left(\id^{\otimes i} \otimes \beta_2 \otimes \id^{\otimes j}\right),
 \end{equation}
where the operators 
\begin{eqnarray*}
&&\beta'_1: \s \RHom(\int[E_n],\int[E_0]) \rightarrow \s \RHom(\int[E_n],\int[E_0]),\\
&&\beta'_2: \s \RHom(\int[E_i],\int[E_{0}])\otimes \s \RHom(\int[E_{n}],\int[E_i])\rightarrow \s \RHom(\int[E_n],\int[E_0]), \\
&&\beta_1: \s \RHom(E_{i+1},E_i) \rightarrow \s \RHom(E_{i+1},E_i) \quad \textrm{and}\\
&&\beta_2: \s \RHom(E_{i+1},E_{i})\otimes \s \RHom(E_{i+2},E_{i+1})\rightarrow \s \RHom(E_{i+2},E_{i}), 
\end{eqnarray*}
are the maps corresponding to the composition and the differential on the spaces of morphisms at the level of the suspension.

It suffices to show that the right and left hand side of \eqref{equationtoprove} evaluate to the same linear
map on all simplices in $ \Pi_{\infty}(A)$. Since the construction is natural with respect to pull-back, we may assume that $A=T\Delta_k$ and that the simplex 
$\sigma$ is the fundamental cochain $[\Delta_k]=(\id: \Delta_k\to \Delta_k)$. Moreover,
by Lemma \ref{lemma:invariance2} we can assume that all the graded vector bundles
$E_i$ are trivial, i.e. $E_i = \Delta_k \times V_i$.

In this situation, Proposition \ref{lemmaMC1} and Proposition \ref{lemmaMC2} provide natural identifications of graded vector spaces:
\begin{eqnarray*}
&&\RHom(E_i,E_j)\cong  \Hom(V_i,V_j)\otimes \Omega(\Delta_k),\\
&&\RHom(\int[E_i],\int[E_j])\cong  \Hom(V_i,V_j)\otimes  C(\Delta_k).
\end{eqnarray*}
Under these identifications, composition of operators corresponds to multiplication in the algebras $\End V\otimes \Omega(\Delta_k)$ and $\End V\otimes \overline{C}(\Delta_k)$, respectively.
Also, the differentials \[\delta:\RHom(E_i,E_j)\rightarrow \RHom(E_i,E_j)\] and 
\[\D:\RHom(\int[E_i],\int[E_j])\rightarrow \RHom(\int[E_i],\int[E_j])\]
correspond to
\begin{eqnarray*}
\Hom(V_i,V_j)\otimes \Omega(\Delta_k)s&\rightarrow &\Hom(V_i,V_j) \otimes \Omega(\Delta_k),\\
\eta &\mapsto& d\eta + \omega^j \wedge \eta -(-1)^{|\eta|}\eta \wedge \omega^i
\end{eqnarray*}
and
\begin{eqnarray*}
\Hom(V_i,V_j) \otimes \overline{C}(\Delta_k) &\rightarrow & \Hom(V_i,V_j)\otimes \overline{C}(\Delta_k),\\
\mu &\mapsto& \overline{\delta} \mu + \alpha^j\overline{\cup} \mu -(-1)^{|\mu|}\mu \overline{\cup}\alpha^i,
\end{eqnarray*}
where $\omega^i$ and $\alpha^i$ are the Maurer-Cartan elements corresponding to the representations up
to homotopy $E_i$ and $\int[E_i]$.

We observe that the map $\int_n$ corresponds to
 
 \[\overline{\int_{n}} : \s (\Hom(V_1,V_0 ) \otimes \Omega(\Delta_k))\otimes \dots \otimes \s (\Hom(V_n,V_{n-1})\otimes \Omega(\Delta_k)) \rightarrow \s (\Hom (V_n, V_0)\otimes \overline{C}(\Delta_k)) \]
 given by the formula
 
 \[(\s \eta^1 \otimes \dots \otimes \s \eta^n)\mapsto  \sum_{l_0\ge 0, \dots, l_n\ge 0}\overline{\psi}^{V}_{l_0+\cdots + l_n + n}
((\s \omega^0)^{\otimes l_0}\otimes \s \eta^1 \otimes \dots \otimes \s \eta^n \otimes (\s \omega^n)^{\otimes l_n}). \]
Let us define $\omega=\omega^0+\dots +\omega^n \in \End V\otimes \Omega(\Delta_k)$ and $\alpha=\alpha_0+\dots +\alpha_n \in \End V\otimes \overline{C}(\Delta_k)$. Clearly,
$\omega $ and $\alpha$ are Maurer-Cartan elements and are related by $\overline{\psi}^{V}(\omega)=\alpha$.

Let $b_2, b'_2$ are the operators corresponding to the products on $\End V\otimes \Omega(\Delta_k)$ and $\End V\otimes \overline{C}(\Delta_k)$ at the level of the suspension
and $b_1, b'_1$ correspond to the twisted differentials $d_\omega$ and $\overline{\delta}_\alpha$.
We need to prove
\begin{equation}\label{toprove}
 b'_1\circ \overline{\int_{n}}+\sum_{i+j=n}b'_2\circ \left( \overline{\int_i} \otimes \overline{\int_j} \right)=\sum_{i+j+1=n}\overline{\int_{n}}\circ \left(\id^{\otimes i} \otimes b_1\otimes \id^{\otimes j}\right)+ \sum_{i+j+2=n}\overline{ \int}_{n-1}\circ \left(\id^{\otimes i} \otimes b_2 \otimes \id^{\otimes j}\right).
 \end{equation}
For this we note that this is exactly the component in $\Hom(V_n,V_0)$ of the equations for the twist $\overline{\psi}^{V}_{\omega}$ of $\overline{\psi} ^{V}$ by $\omega$
to be an $\mathsf{A}_{\infty}$ morphism.
\end{proof}

\section{Examples}\label{section:examples}

In the following we describe some examples of representations up to homotopy
of Lie algebroids to which the $\A_\infty$ functor
\begin{align*}
\int: \RRep(A) \to \URRep(\Pi_{\infty}(A))
\end{align*}
can be applied.

\subsection{Parallel transport for superconnections}

The situation where the Lie algebroid is the tangent bundle $TM$ of a manifold $M$  has been studied 
by Igusa $\cite{I}$ and Block-Smith \cite{BS}, and was the starting point of  our work.
In this case, representations up to homotopy are precisely the
$\mathbb{Z}$-graded versions of Quillen's flat superconnections (\cite{Quillen2}).
Since  \[\mathrm{Hom}_{\mathsf{DGCA}}(\Omega(TM),\Omega(\Delta_k))\cong\mathrm{Hom}_{\smooth}(\Delta_k, M),\] 
we know that

\[\Pi_{\infty}(TM) \cong \mathrm{Sing}(M).\]

Thus, integrating a representation up to homotopy of the tangent Lie algebroid $TM$ amounts to assigning holonomies to smooth singular chains on $M$.
As explained in \cite{BS}, this procedure is a generalization of the Riemann-Hilbert correspondence. Just as flat connections correspond, via holonomy, to representations
of the fundamental groupoid of $M$, flat superconnections correspond to representations up to homotopy of the infinity groupoid of $M$.

\subsection{Ordinary representations}\label{sub:compatibility}

Let $A$ be an integrable Lie algebroid with integrating source simply connected Lie groupoid
$G$. By Lie's second theorem for Lie algebroids, every one-simplex
\[\sigma \in \mathrm{Hom}_{\mathsf{Lie-alg}}(TI, A)=\mathrm{Hom}_{\mathsf{DGCA}}(\Omega(A), \Omega(I))\]
integrates to a morphism of Lie groupoids
\[\tilde{\sigma}: I \times I \rightarrow G\]
from the pair groupoid of $I$ to $G$. We set the target (source) map to be the projection
onto the first (second) factor, respectively.

Let $\tau(\sigma)$ be the value of $\tilde{\sigma}$ on $(1,0) \in I\times I$, i.e.
\begin{align*}
\tau(\sigma):=\tilde{\sigma}(1,0).
\end{align*}
One obtains a morphism of simplicial sets
\[\pi: \Pi_{\infty}(A) \rightarrow \mathsf{N} G\]
from the $\infty$-groupoid $\Pi_{\infty}(A)$ of $A$ to the nerve $\mathsf{N} G$ of the Lie groupoid $G$
by setting
\[\pi(\sigma):=(\tau(\sigma_1),\dots,\tau(\sigma_k)).\]
Here, $\sigma$ is a $k$-simplex in $\Pi_{\infty}(A)$ and $\sigma_i$ is the one simplex given by the edge joining $v_{i}$ and $v_{i-1}$; that is: the one simplex given by pre-composing $\sigma$ with
(the tangent map of)
\begin{align*}
t \mapsto (\underbrace{1,\dots,1}_{i-1\text{ times}},1- t, 0,\dots,0).
\end{align*}

The morphism $\pi: \Pi_{\infty}(A) \to \mathsf{N}G$ induces a functor of $dg$-categories
\begin{align*}
\pi^{*}: \URRep(G) \hookrightarrow \URRep(\mathsf{N} G) \to \URRep(\Pi_{\infty}(A)).
\end{align*}
The inclusion $\URRep(G) \hookrightarrow \URRep(\mathsf{N}G)$ reflects the fact
that we require the graded vector spaces $(E_x)_{x\in M}$ underlying a representation up to homotopy of a {\em Lie} groupoid
$G$ over $M$ to fit into a graded {\em vector bundle} $E$.

\begin{proposition}\label{proposition:ordinary_reps}
Let $E$ be an ordinary representation of the Lie algebroid $A$ and denote by $\tilde{E}$ the corresponding representation of $G$.

Then
\[\pi^*\tilde{E}=\int [E]\]
holds
in $\URRep(\Pi_{\infty}(A))$.
\end{proposition}

\begin{proof}
By degree reasons we know that the only nontrivial structure operator for $\int [E]$ is $F_1$.
Clearly, the same is true for $\pi^*(\tilde{E})$.
Let us denote by $\mu:G \rightarrow \End E$ the structure operator for the representation $\tilde{E}$. We need to prove that for any one simplex
$\sigma:TI \rightarrow A$ one has
\[F_1(\sigma)=\mu(\pi(\sigma)).\]

Naturality of the the integration construction -- see Remark \ref{remark:natural} --
implies
\[F_1(\sigma)=F_1^E(\sigma \circ \id)= F^{\sigma^*(E)}_1(\id).\]
We denote by $\theta: I \times I \rightarrow \End(\sigma^*E)$ the structure map corresponding to the representation that integrates $\sigma^*E$ and observe that since
$\theta= \mu \circ \tilde{\sigma}$ we know that
\[\theta(1,0) = \mu(\tilde{\sigma}(1,0)) = \mu(\widetilde{\sigma_1}(1,0)) = \mu(\pi(\sigma))\]
holds.
Thus, it is sufficient to prove that
 \[F^{\sigma^*(E)}_1(\id)=\theta(1,0),\]
i.e. we may assume without loss of generality $A=TI$ and $\sigma=\id$.
In this case the statement reduces to the fact that
the holonomy of a path is given by parallel transport along the direction-reversed path, see Remark \ref{remark:parallel_transport}. 
\end{proof}

\begin{remark}
The differentiation functor
\[ \Psi:\URRep(G) \rightarrow \RRep(A)\]
constructed in  \cite{AS} can be extended naturally to a $dg$-functor.
The definition of $\Psi$ involves iterated differentiation and in a sense it is given by the inverse construction of the integration procedure described above.
Proposition \ref{proposition:ordinary_reps} shows that, when restricted to ordinary representations, this is indeed the case. However, we would like to point out that the diagram
\begin{eqnarray*}
\xymatrix{
 \URRep(G) \ar[rr]^(0.45){\pi^*}\ar[dd]_{\Psi} && \URRep(\Pi_{\infty}(A))\\
&&\\
 \RRep(A)\ar[rruu]_{\int}& & 
}
\end{eqnarray*}
does not commute in general. We expect that it commutes up to homotopy, but this issue will not be discussed here.
\end{remark}

\subsection{Lie algebras and flat connections}

Let us fix a finite dimensional Lie algebra $\frak{g}$ and denote by
$G$ the simply connected Lie group integrating it. We think of the Chevalley-Eilenberg complex $CE(\frak{g})$ as the
space of functions on the $dg$-manifold $\frak{g}[1]$ determined by $\frak{g}$. 

The space of $k$-simplices on $\frak{g}$
\[\Pi_{\infty}(\frak{g})_k= \mathrm{Hom}_{\mathsf{DGCA}}(CE(\frak{g}),\Omega(\Delta_k))\]
can be naturally identified with the space of flat connections on the trivial $G$ bundle over $\Delta_k$. In fact, any such map is in particular
a map of graded commutative algebras, and since $CE(\frak{g})$ is free as a graded commutative algebra, one obtains
\[ \mathrm{Hom}_{\mathsf{GCA}}(CE(\frak{g}),\Omega(\Delta_k))\cong  \mathrm{Hom}_{\mathsf{Vect}}(\frak{g^*},\Omega^1(\Delta_k)) \cong \frak{g}\otimes \Omega^1(\Delta_k).\]
One easily checks that the condition of commuting with the differential corresponds to the Maurer-Cartan equation.
We conclude that integrating a representation up to homotopy of a Lie algebra amounts to assigning holonomies to the spaces of flat connections on the trivial $G$-bundles over the simplices.

\subsection{The adjoint representation of a Lie algebroid}

Given a Lie algebroid $A$, there is a representation up to homotopy $\mathrm{ad}(A) \in \RRep(A)$ -- called the adjoint representation of $A$ -- which is well defined up to isomorphism. 
For instance:
\begin{itemize}
\item If $A$ is a Lie algebra, $\mathrm{ad}(A)$ coincides with the usual adjoint
representation.
\item If $A$ is a foliation, $\mathrm{ad}(A)$ is quasi-isomorphic to the Bott connection.
\end{itemize}
The cohomology associated to the adjoint
representation measures the deformations of the Lie algebroid structure, as one would expect from the Lie algebra case -- see \cite{ AC1,AS,CrM} for more details.

In order to define the adjoint representation, one has to choose a connection $\nabla$ on the vector bundle $A$. The adjoint representation of $A$ induced
by $\nabla$ is the representation up to homotopy of $A$ on the vector bundle $A \oplus TM$, where $A$
is in degree zero and $TM$ in degree one, given by the following structure operator:
 \begin{equation*}
D=\partial+\nabla^{bas}+K_{\nabla}.
\end{equation*}
 The differential $\partial$ on the graded vector bundle is the anchor map of $A$, $\nabla^{bas}$ is an $A$-connection and $K_{\nabla}$ is an
 endomorphism-valued cochain. The latter two are defined by the formulae
 \begin{eqnarray*}
 \nabla^{bas}_{\alpha}(\beta)&:=& \nabla_{\rho(\beta)}(\alpha)+ [\alpha, \beta],\\
 \nabla^{bas}_{\alpha}(X)&:=& \rho(\nabla_{X}(\alpha))+ [\rho(\alpha), X] \quad ,\\
 K_{\nabla}(\alpha,\beta)(X)&:=&\nabla_X([\alpha,\beta])-[\nabla_X(\alpha),\beta]-
[\alpha,\nabla_X(\beta)]-\nabla_{\nabla_{\beta}^{\textrm{bas}}X}(\alpha)
+\nabla_{\nabla_{\alpha}^{\textrm{bas}}X}(\beta).
\end{eqnarray*}
This representation up to homotopy is denoted by $\ad_{\nabla}(A)$. The isomorphism class of this
representation, which is independent of the connection $\nabla$, is denoted by $\ad(A)$.

\subsection{Poisson manifolds}

Let $P$ be a Poisson manifold with Poisson bivector field $\pi \in \Gamma(\wedge^2 TP)$.
The cotangent bundle $T^*P$ has the structure of a Lie algebroid with anchor map $TP^*\rightarrow TP$
given by contraction with $\pi$ and bracket determined by the formula
\[[df,dg]=d\{f,g\}.\]

The Poisson sigma model (\cite{Ikeda,Strobl, CF}) is a two dimensional topological field theory associated to a Poisson manifold, whose perturbative quantization
gives Kontsevich's formula for a $\star$-product of $\pi$.
Let us briefly recall the main ingredients of this theory:
The  fields on an oriented surface $\Sigma$ are the vector
bundle maps $\hat{X}: T\Sigma \rightarrow T^*P$.
It is convenient to write $\hat{X}=(X, \eta)$, where $X: \Sigma \to P$ is the base map and $\eta \in \Gamma(\mathrm{Hom}(T\Sigma, X^*(T^*P)))$.
The action functional is

\[ S(X,\eta):=\int_{\Sigma} \langle \eta,dX\rangle +\frac{1}{2}\langle \pi \circ X,\eta \wedge \eta \rangle. \]
The classical solutions of this theory are those vector bundle maps which are Lie algebroid morphisms from $T\Sigma$ to $T^{*}P$.

We conclude that integrating a representation up to homotopy of the Lie algebroid $T^*P$ amounts to assigning holonomies
to the classical fields of the Poisson sigma model on the simplex.
We hope to make the connection between the integration of representations
up to homotopy of $T^{*}P$ and the Poisson sigma model more precise elsewhere.

\subsection{$\Pi_{\infty}(-)$ vs $\Pi_1(-)$}\label{subsection:Pi_infty}

\begin{remark}
In Section \ref{section:functor}, the integration $\int[E]$
of a representation up to homotopy $E$ of a Lie algebroid $A$ was defined. It is a representation up to homotopy
of the $\infty$-groupoid $\Pi_{\infty}(A)$ of $A$.

Let us assume that $A$ can be integrated to a source simply connected Lie groupoid $G$, which is Hausdorff.
In Subsection \ref{sub:compatibility}, a map of simplicial sets
\begin{align*}
\pi: \Pi_{\infty}(A) \to \mathsf{N}G
\end{align*}
was constructed.
One might wonder whether the above integrability assumption on $A$ implies that $\int$ factors through $\pi^{*}$, i.e. whether it is possible to complete the diagram
\begin{align*}
\xymatrix{
 \URRep(G) \ar[rr]^{\pi^*\,\,\,\,} && \URRep(\Pi_{\infty}(A))\\
&& \\
 \RRep(A)\ar[rruu]_{\int} \ar@{-->}[uu]^{\exists ?}&&\\
}
\end{align*}
in a commutative way. Intuitively, one should not expect this to be the case, since $G$ only contains information about homotopy classes of $A$-paths and not about the 
higher homotopies. In any case, we present an example that shows that this factorization is not possible.
\end{remark}

\begin{proposition}\label{proposition:dont_factor}
The integration $\int[E]$ of the representation up to homotopy $E$ of $TS^{2}$ associated to
 the graded vector bundle $E = E^{0}\oplus E^{1}$ given by the trivial vector bundle $S^{2}\times \mathbb{R}$ in both degrees,
  equipped with the structure operator
\begin{align*}
D = d + \eta,
\end{align*}
where $\eta \in \Omega^{2}(S^{2})$ is a volume form, viewed as an element of
$\Omega^{2}(S^{2},\Hom(E^{1},E^{0})),$
is not quasi-isomorphic to any element of $\pi^{*}\URRep(G)$.
\end{proposition}

\begin{remark}
In order to prove the proposition, we need to introduce the cohomology associated to a representation
up to homotopy $E$ of a simplicial set $X_{\bullet}$:
\end{remark}

\begin{definition}
Let $E$ be a representation up to homotopy of $X_{\bullet}$.
The  cohomology of $X_{\bullet}$ with values in $E$ is the graded vector space
\begin{align*}
H(X_{\bullet},E) := H(\RHom(\mathbb{R},E)), 
\end{align*}
where $\mathbb{R}$ denotes the trivial representation up to homotopy of $X_{\bullet}$ on $\mathbb{R}$.
\end{definition}

We can now proceed with the proof of  Proposition \ref{proposition:dont_factor}.

\begin{proof}
Observe that the source simply connected Lie groupoid integrating $TS^{2}$ is the pair groupoid
$S^{2}\times S^{2}$. 
Next, suppose that $\mathcal{E}$ is a representation up to homotopy of $G$, whose pull back
$\pi^{*}\mathcal{E}$ is quasi-isomorphic to $\int[E]$. By definition, this implies that the fiberwise cohomologies
are isomorphic, hence
\begin{align*}
{\cal H}(\mathcal{E}_x)= {\cal H}((\pi^*\mathcal{E})_x)\cong {\cal H}(\int[E]_x)= {\cal H}(E_x)= E_x=\mathbb{R}\oplus \mathbb{R}[-1].
\end{align*}
Since any representation up to homotopy can be transfered to its cohomology vector bundle, we may assume
without loss of generality that $\mathcal{E}_x = \mathbb{R}\oplus \mathbb{R}[-1]$.

By degree reasons we know that the structure operator of $\mathcal{E}$ is of the form
\[ F_1+ F_2.\]
Next,  we claim that one can choose an automorphism of the vector bundle $\mathcal{E}$ that conjugates $F_1$ to the identity. Indeed, choose a point $x \in S^2$ and define the automorphism
\[\phi: \mathcal{E} \rightarrow \mathcal{E}\] 
by setting
\[\phi(y)=F_1(y,x).\]
Since $V$ is unital, one can easily check that $\id = \phi^{-1} \circ F_1\circ \phi$. Thus, we can further assume $F_1=\id$.

The operator $F_2 \in C^2(G)$ is a closed two-cocycle of the pair groupoid $S^2 \times S^2$. Since the nerve of this groupoid is contractible, we know that $F_2$ is
exact and therefore $\mathcal{E}$ is isomorphic to the sum of trivial representations $\mathbb{R} \oplus \mathbb{R}[-1]$. This implies
\[\pi^*(\mathcal{E}) \cong \underline{\mathbb{R}} \oplus \underline{\mathbb{R}}[-1].\]

We can now compute the cohomology associated to this representation:
\[H(\Pi_{\infty}(TS^2), \pi^*\mathcal{E})\cong H(\textrm{Sing}(S^2), \mathbb{R} \oplus \mathbb{R}[-1])\cong H(S^2)\oplus H(S^2)[-1].\]
In particular, we observe that 
\[H^2(\Pi_{\infty}(TS^2), \pi^*\mathcal{E})\cong \mathbb{R}.\]
Since quasi-isomorphisms induce isomorphisms in cohomology, it is enough to prove that
 \[H^2(\Pi_{\infty}(TS^2), \int[E])\]
vanishes,
in order to contradict our assumption that $\pi^{*}\mathcal{E}$ and $\int[E]$ are quasi-isomorphic.

Simple degree considerations imply that the only nontrivial structure operator of $\int[E]$ is 
$F'_2\in C^2(S^2)$. We claim that $[F'_2]\neq 0 \in H^2(S^2)$. Indeed:
\[F'_2(\sigma)=\pm \int_I\Theta^*_{(2)}(\Path \sigma)^{*}(\Chen(\s \eta))=\pm \int_{\Delta_2}\sigma^{*}\eta.\]
By de Rham's theorem, we conclude that $F'_2$ is a nontrivial cohomology class.

To compute the cohomology of $\mathrm{Sing}(S^2)$ with coefficients in $\int [E]$, we consider the filtration of $C(\mathrm{Sing}(S^2), \int[E])$ given
by the cochain degree. It induces a spectral sequence which converges to $H(\textrm{Sing}(S^{2}),\int [E])$ and whose second page is

\begin{equation*}
\frame{ \xymatrix{
0&0&\,\,\,\,\,\,0\\
H^2(S^2)&H^2(S^2)& \, \,\,\,\,\,0\\
0& 0&\, \,\, \,\,\,0\\
H^0(S^2)&H^0(S^2)\ar[luu]_{d_2}&\,\,\,\,\,\,0
}}
\end{equation*}

The operator $d_2$ is given by multiplication with $[F'_2]$ and since this class is not zero, the second cohomology
\[H^2(\Pi_{\infty}(TS^2),\int [E])\]
vanishes.
\end{proof}

\begin{remark}
Mimicking the arguments from above, one can prove that there is no representation up to homotopy $\mathcal{E}$
of $G$, whose image under the differentiation map
\begin{align*}
\Psi: \URRep(G) \to \RRep(A)
\end{align*}
from \cite{AS} equals $E$.
In fact, if such an $\mathcal{E}$ would exist, the cohomology class
of $\omega$ would be in the image of the van Est map and would hence vanish,
since $H^{2}(S^{2}\times S^{2}) = 0$.

\end{remark}

\appendix

\section{Appendix}\label{appendix}

\subsection{Representations up to homotopy of Lie algebroids }\label{app1}

We review the definitions and basic facts regarding
representations up to homotopy of Lie algebroids. More details on these constructions, as well as the proofs of the
results stated here, can be found in
\cite{AC1,AS}. Throughout the appendix $A$ denotes a Lie algebroid over a manifold
$M$.

\begin{remark}
Given a Lie algebroid $A$, there is a differential graded algebra
$\Omega(A)= \Gamma(\wedge A^*)$, with differential defined via
the Koszul formula
\begin{eqnarray*}
d \omega(\alpha_1, \dots ,\alpha_{n+1})&=& \sum_{i<j}(-1)^{i+j} \omega([\alpha_i,\alpha_j],\cdots,\hat{\alpha}_i,\dots,\hat{\alpha}_j,\dots,\alpha_{k+1})\\
&&+\sum_i(-1)^{i+1}L_{\rho(\alpha_i)}\omega(\alpha_1,\dots,\hat{\alpha}_i,\dots,\alpha_{k+1}),
\end{eqnarray*}
where $\rho$ denotes the anchor map and $L_{X}(f)= X(f)$ is the Lie derivative along vector fields.
The operator $d$ is a coboundary operator ($d^{2}= 0$) and satisfies
the derivation rule
\[ d(\omega\eta)= d(\omega)\eta+ (-1)^p\omega d(\eta),\]
for all $\omega\in \Omega^p(A), \eta\in \Omega^q(A)$.

Given a graded vector bundle $E=\bigoplus_{k\in \mathbb{Z}}E^k $ over $M$, we denote by $\Omega(A,E)$ the space 
\begin{align*}
\Gamma(E\otimes \wedge A^*),
\end{align*}
graded with respect to the total degree. The wedge product gives this space the structure of a graded commutative module over the algebra $\Omega(A)$, i.e.
$\Omega(A,E)$ is a bimodule over $\Omega(A)$ and
\begin{equation*}
\omega \wedge \eta=(-1)^{kp} \eta \wedge \omega
\end{equation*}
holds for $\omega \in \Omega^k(A)$ and $\eta \in \Omega(A,E)^p$. 
For a detailed explanation of the sign conventions used in the definition in the wedge product, see the appendix of \cite{AC1}. 
In order to simplify the notation we will sometimes omit the wedge symbol.
\end{remark}

\begin{definition} A representation up to homotopy of $A$ consists of a graded vector bundle $E$ over $M$
and a linear operator
\begin{equation*}
D:\Omega(A,E)\rightarrow \Omega(A,E)
\end{equation*}
which increases the total degree by one and satisfies $D^2= 0$, as well as
the graded derivation rule
\begin{equation*}
D(\omega \eta)=d(\omega) \eta +(-1)^k \omega D(\eta)
\end{equation*}
for all $\omega\in \Omega^k(A)$, $\eta\in \Omega(A, E)$.
\end{definition}

\begin{remark}
The representations up to homotopy of  $A$ can be naturally organized into 
a $dg$-category $\RRep(A)$: A morphism $\phi: E\to F$ between two representations up to homotopy
of $A$ is a map
\begin{equation*}
\phi:\Omega(A,E) \rightarrow \Omega(A,F)
\end{equation*}
which is $\Omega(A)$-linear in the graded sense. We denote the space of 
all morphisms from $E$ to $F$ by $\RHom(E,F)$. This is a graded vector space and we denote
by $\RHom^k(E,F)$ the subspace of homogeneous elements of degree $k$. There is a differential
\begin{equation*}
\delta:  \RHom^k(E,F) \rightarrow  \RHom^{k+1}(E,F)
\end{equation*}
given by the formula 
\begin{equation*}
\delta(\phi)=D_F \circ \phi -(-1)^k \phi \circ D_E.
\end{equation*}

Let $\mathrm{Hom(E,F)}$ be 
the space of morphisms of degree zero which
commute with the differentials. Observe that
\begin{equation*}
\mathrm{Hom}(E,F)=Z^0(\RHom(E,F)).
\end{equation*}
\end{remark}

\begin{remark}
Let $E$ be a graded vector space. The space $\Omega(A,E)= E \otimes \Omega(A)$ is a differential left  module over the $dg$-algebra
$\End E\otimes \Omega(A)$. Given $\omega \in \End E\otimes \Omega(A)$, we denote the corresponding
operator on $E \otimes \Omega(A)$ by $\omega \wedge$. 
\end{remark}
One can describe representations up to homotopy  on trivial vector bundles as follows:

\begin{proposition}\label{lemmaMC1}
Let $A$ be a Lie algebroid over $M$ and $V$ a finite dimensional graded vector space. 
\begin{itemize}
\item[a)]
There is a natural bijective correspondence
between
\begin{enumerate}
\item Maurer-Cartan elements of the $dg$-Lie algebra $\End E\otimes \Omega(A)$.
\item Representations up to homotopy of $A$ on the trivial vector bundle $M\times V$.
\end{enumerate}
The correspondence is given by
\[\omega \mapsto d+\omega\wedge \,.\]

\item[b)]
Suppose that $E$ and $E'$ are representations up to homotopy on trivial vector bundles $M\times V$ and $M\times V'$, respectively. There is a natural isomorphism
\begin{eqnarray*}
 \Hom(V,V')\otimes \Omega(A)&\cong& \RHom(E,E').\\
 \eta &\mapsto& \eta\wedge 
\end{eqnarray*}
Under this identification, the operator $\delta:\RHom(E,E')\rightarrow \RHom(E,E')$ corresponds to the map
\begin{eqnarray*}
 \Hom(V,V')\otimes \Omega(A)&\rightarrow &  \Hom(V,V')\otimes \Omega(A) \\
 \eta &\mapsto& d \eta +\omega' \wedge \eta -(-1)^{|\eta|}\eta \wedge \omega,
\end{eqnarray*}
where $\omega$ and $\omega'$ are the Maurer-Cartan elements corresponding to $E$ and $E'$ and the wedge product 
is taken in the algebra $\End(V\oplus V')\otimes \Omega(A)$.
Furthermore, composition corresponds to multiplication in the algebra $\End(V\oplus V'\oplus V'')\otimes \Omega(A)$.
\end{itemize}
\end{proposition}

\subsection{$\mathsf{A}_{\infty}$-algebras and morphisms}\label{app2}

The notion of  an $\mathsf{A}_{\infty}$-algebra was introduced in the sixties by J. Stasheff \cite{Sta1, Sta2}, and has since proved to be important in several areas of mathematics. 
In order to fix our conventions, we collect here some definitions and basic facts regarding $\mathsf{A}_{\infty}$-algebras. 

\begin{definition}
Let $V=\bigoplus_{k \in \mathbb{Z}}V^k$ be a graded vector space. The suspension $\mathsf{s}V$ of $V$ is the graded vector space given by
\[ (\mathsf{s}V)^k := V^{k+1}.\]
\end{definition}

\begin{definition}
An $\mathsf{A}_{\infty}$-algebra is a graded vector space $\mathsf{A}$, together with a sequence of linear maps of degree one
\[b_n: (\mathsf{s}\mathsf{A})^{\otimes n}\rightarrow s\mathsf{A}, \]
satisfying, for each $n\geq 1$, the equations:
\begin{equation}
 \sum_{i+j+k=n} b_{i+k+1} \circ (\id^{\otimes i}\otimes b_j \otimes \id^{\otimes k})  =0. 
 \end{equation}
\end{definition}

\begin{examples}
If the graded vector space $\mathsf{A}$ is concentrated in degree zero, an $\mathsf{A}_{\infty}$-algebra on $\mathsf{A}$ is the same as an associative algebra.
A differential graded algebra is the same as an $\mathsf{A}_{\infty}$-algebra where $b_n=0$ for $n\notin\{1,2\}$.

Observe that $b_1$ always defines a coboundary operator on the graded vector space $\mathsf{A}$.
\end{examples}

\begin{remark}
An alternative definition of $\mathsf{A}_{\infty}$-algebras can be described in terms of maps
$\tilde{b}_n : \mathsf{A}^{\otimes n} \rightarrow \mathsf{A}. $
The advantage of using the suspension is that no signs appear in the structure equations. The definition in terms of the 
structure maps $\tilde{b}_n$ requires some signs that can be determined by requiring the following diagram to commute:

\begin{equation*}
\xymatrix{
(\mathsf{sA})^{\otimes n} \ar[r]^{b_n} \ar[d]_{(s^{-1})^{\otimes}}&\mathsf{sA}\\
\mathsf{A}^{\otimes n}\ar[r]^{\tilde{b}_n}&\mathsf{A}. \ar[u]_s\\
}
\end{equation*}

Here, the map $\mathsf{s}:\mathsf{A} \rightarrow \mathsf{sA}$ is given by $\mathsf{s}(v)=v\in \mathsf{sA}$.
\end{remark}

\begin{definition}
Let $\mathsf{A}$ and $\mathsf{A}'$ be two $\mathsf{A}_{\infty}$-algebras. A morphism $\psi:\mathsf{A} \rightarrow \mathsf{A}'$ is a sequence of degree zero maps

\[ \psi_n:(\mathsf{sA})^{\otimes n}\rightarrow \mathsf{sA},  \] 
satisfying, for each $n \geq 1$,  the equation
\begin{equation}\label{morphismainfinity}
 \sum_{i+j+k=n} \psi_{i+k+1}\circ (\id^{\otimes i}\otimes b_j\otimes \id^{\otimes k}  )= \sum_{l_1+\dots+l_r=n} b'_{r}\circ (\psi_{l_1}\otimes \dots \otimes \psi_{l_r}  ). 
\end{equation}
A morphism $\psi: \mathsf{A} \to \mathsf{A}'$ is an $\mathsf{A}_{\infty}$ quasi-isomorphism if the chain map
\begin{align*}
\psi_1: (\mathsf{A}, b_1) \to (\mathsf{A}',b'_1)
\end{align*}
induces an isomorphism in cohomology.

If $\psi: \mathsf{A} \rightarrow \mathsf{A}'$ and $\psi':\mathsf{A}' \rightarrow \mathsf{A}''$ are to morphisms, their composition $\psi' \circ \psi$ is defined by
\[(\psi' \circ \psi)_n= \sum_{i_1+\dots +i_r=n} \psi'_r\circ (\psi_{i_1} \otimes \dots \otimes \psi_{i_r}).\]
The identity morphism $ \id:\mathsf{A} \rightarrow \mathsf{A}$ is defined by setting $\id_1=\id$ and $\id_n=0$ for $n\neq 1$.
\end{definition}

\begin{remark}
In the case where $\psi: \mathsf{A} \rightarrow \mathsf{A}'$ is an $\mathsf{A}_{\infty}$ morphism between $dg$-algebras, the structure equations (\ref{morphismainfinity}) take the form:

\begin{equation}\label{equations for functor}
 b'_1\circ \psi_{n}+\sum_{i+j=n}b'_2\circ (\psi_i \otimes \psi_j)=\sum_{i+j+1=n}\psi_n\circ (\id^{\otimes i} \otimes b_1\otimes \id^{\otimes j})+ \sum_{i+j+2=n} \psi_{n-1}(\id^{\otimes i} \otimes b_2 \otimes \id^{\otimes j})
\end{equation}
\end{remark}

\begin{definition}
Let $\mathsf{A}$ be an $\mathsf{A}_{\infty}$-algebra such that $b_n=0$, except for finitely many values of $n$. 
A Maurer-Cartan element of $\mathsf{A}$ is an element $ x \in (\mathsf{sA})^0$ satisfying
\[\sum_{n\geq1}b_n(x^{\otimes n})=0.\]
We will sometimes abuse our conventions and say that $y \in \A^1$ is a Maurer-Cartan element if $\s y\in (\s \A)^0$ is a Maurer-Cartan element.
\end{definition}

\begin{proposition}\label{MCgoestoMC}
Let $\psi:\mathsf{A}\rightarrow \mathsf{A}'$ be an $\mathsf{A}_{\infty}$ morphism between $\mathsf{A}_\infty$-algebras such that $b_n$, $b'_n$ and $\psi_n$ are nonzero
only for finitely many values of $n$.
If $x \in (\mathsf{sA})^0$ is a Maurer-Cartan element of $\mathsf{A}$, then

\[\psi(x) := \sum_{n\geq 1}\psi_{n}(x^{\otimes n})\]
 is a Maurer-Cartan element of $\mathsf{A}'$.
\end{proposition}

\begin{remark}
The hypothesis that only finitely many of the structure operators are nonzero guarantees the convergence of the infinite sum above. 
The statement remains true for
arbitrary $\mathsf{A}_{\infty}$-algebras and morphisms, provided that the infinite sums converge appropriately. Similar issues arise in the definition of twisting by a Maurer-Cartan element given below.
\end{remark}

\begin{definition}
Let $\mathsf{A}$ be an $\mathsf{A}_{\infty}$-algebra and $x \in (\s \mathsf{A})^{0}$ a Maurer-Cartan element.

The $\mathsf{A}_{\infty}$-algebra $\mathsf{A}_x$ is given by the structure maps
\begin{align*}
(b_{x})_n(a_1\otimes \cdots \otimes a_n) := \sum_{l_0\ge 0, \cdots ,l_n\ge 0} b_{n+l_0+\cdots + l_n}(x^{\otimes l_0} \otimes a_1 \otimes x^{\otimes l_1}\otimes \cdots \otimes x^{\otimes l_{n-1}} \otimes a_n \otimes x^{\otimes l_n}).
\end{align*}
Here, $\{b_n\}$ denotes the family of structure maps of the original $\mathsf{A}_{\infty}$-algebra $\mathsf{A}$.

We say that $\mathsf{A}_x$ is obtained from $\mathsf{A}$ by twisting by $x$.
\end{definition}

\begin{definition}
Suppose $\psi: \mathsf{A} \to \mathsf{A}'$ be a morphism of $\mathsf{A}_{\infty}$-algebras and let $x \in (\s \mathsf{A})^{0}$ be a Maurer-Cartan element.
There is a morphism of $\mathsf{A}_{\infty}$-algebras
\begin{align*}
\psi_x: \mathsf{A}_x \to \mathsf{A}'_{\psi(x)}
\end{align*}
between the twisted $\mathsf{A}_{\infty}$-structures, given by
\begin{align*}
(\psi_x)_n(a_1\otimes \cdots \otimes a_n) := \sum_{l_0 \ge 0, \cdots l_n \ge 0} \psi_{n+l_0+\cdots +l_n}(x^{\otimes l_0} \otimes a_1 \otimes x^{\otimes l_1}\otimes \cdots \otimes x^{\otimes l_{n-1}} \otimes a_n \otimes x^{\otimes l_n}).
\end{align*}
We say that $\psi_x$ is obtained from $\psi$ by twisting with $x$.
\end{definition}

Constructing tensor products of $\mathsf{A}_\infty$-algebras and morphisms is a complicated issue in general.
For our purposes we will only need the following special case:

\begin{proposition}
Let $\mathsf{A}$ and $\mathsf{A}'$ be differential graded algebras and $\psi: \mathsf{A}\rightarrow \mathsf{A}'$ an $\mathsf{A}_\infty$ morphism between them. Then, for
any graded algebra $E$, we define a sequence $\phi$ of maps
$\phi_1,\phi_2,\dots, $
where
\[\phi_n: \mathsf{s}(E \otimes\mathsf{A})^{\otimes n}\rightarrow  \mathsf{s}(E \otimes\mathsf{A'})\]
is defined by the formula

\[\phi_n((e_1\otimes \s a_1)\otimes \dots\otimes (e_n\otimes \ s a_n)):=(-1)^{\sum_{i=1}^n [a_i](|e_{i+1}|+ \dots +|e_{n}|)}(e_1 \dots e_n) \otimes \psi_n(\s a_1 \otimes \dots \otimes \s a_n).\]\
Here, we see $e\otimes \s a$ as an element of $\mathsf{s}(E \otimes \mathsf{A})$ via
\begin{align*}
\s (E\otimes A) \cong E \otimes \s \A, \quad \s (e\otimes a) \mapsto (-1)^{|e|} e \otimes \s a.
\end{align*}

The family $\phi$ defines an $\mathsf{A}_\infty$ morphism from $E \otimes \mathsf{A}$ to $E \otimes \mathsf{A'}$, which we denote by $\id_E \otimes \psi$.
\end{proposition}

\begin{remark}
The notion of an $\mathsf{A}_\infty$ algebra generalizes naturally to that of an $\mathsf{A}_{\infty}$-category in such a way that an $\mathsf{A}_\infty$-algebra is an $\mathsf{A}_{\infty}$-category
with only one object. For the purpose of the present paper we will only need to consider $\mathsf{A}_{\infty}$ functors between $dg$-categories:
\end{remark}

\begin{definition}
Let $\mathsf{C}$ and $\mathsf{C'}$ be differential graded categories. An $\mathsf{A}_{\infty}$ functor $\psi: \mathsf{C}\rightarrow \mathsf{C'}$ consists of the following data:
\begin{enumerate}
\item A function $f: Ob(\mathsf{C})\rightarrow Ob(\mathsf{C'})$.
\item For $n\geq 1$ and any sequence of objects $v_0,\dots v_n$, a degree zero map
\[ \psi_n: \mathsf{s}\RHom(v_1,v_{0}) \otimes \dots \otimes  \mathsf{s}\RHom(v_n,v_{n-1})\rightarrow \RHom(f(v_n),f(v_0)),  \]
satisfying equations \eqref{equations for functor}, as well as $\psi_1(\s(\id)) = \s(\id)$
and $\psi_n(\cdots \otimes \s(\id) \otimes \cdots) = 0$ for $n>1$.

\end{enumerate}

\end{definition}

\end{document}